\documentclass[11pt]{article}

\title{Lecture notes on locally compact quantum groups}
\author{Alfons Van Daele}

\usepackage{amsmath}  
\usepackage{amsthm} 
\usepackage{amssymb} 

\usepackage{hyperref} 

\usepackage{color} 

\newtheorem{thm}{Theorem}
\newtheorem{inspr}[thm]{}

\newenvironment{definitie}{\begin{itemize}\item[ ]\hspace{-26pt}\bf Definition \rm }{\end{itemize}}
\newenvironment{notatie}{\begin{itemize}\item[ ]\hspace{-26pt}\bf Notation \rm }{\end{itemize}}
\newenvironment{voorbeeld}{\begin{itemize}\item[ ]\hspace{-26pt}\bf Example \rm }{\end{itemize}}
\newenvironment{stelling}{\begin{itemize}\item[ ]\hspace{-26pt}\bf Theorem \rm }{\end{itemize}}
\newenvironment{propositie}{\begin{itemize}\item[ ]\hspace{-26pt}\bf Proposition \rm }{\end{itemize}}
\newenvironment{lemma}{\begin{itemize}\item[ ]\hspace{-26pt}\bf Lemma \rm }{\end{itemize}}
\newenvironment{opmerking}{\begin{itemize}\item[ ]\hspace{-26pt}\bf Remark \rm }{\end{itemize}}
\newenvironment{voorwaarde}{\begin{itemize}\item[ ]\hspace{-26pt}\bf Condition \rm }{\end{itemize}}

\newenvironment{commentaar-0}{\blauw\it} 

\newenvironment{commentaar-1}{\rood\it}   

\newenvironment{commentaar-2}{\magenta\it} 

\renewcommand{\Bbb}{\mathbb} 

\newcommand{\defin}{\begin{inspr}\begin{definitie}}  
\newcommand{\edefin}{\end{definitie}\end{inspr}}

\newcommand{\notat}{\begin{inspr}\begin{notatie}}  
\newcommand{\enotat}{\end{notatie}\end{inspr}}

\newcommand{\voorb}{\begin{inspr}\begin{voorbeeld}}  
\newcommand{\evoorb}{\end{voorbeeld}\end{inspr}}

\newcommand{\stel}{\begin{inspr}\begin{stelling}}
\newcommand{\estel}{\end{stelling}\end{inspr}}

\newcommand{\prop}{\begin{inspr}\begin{propositie}}
\newcommand{\eprop}{\end{propositie}\end{inspr}}

\newcommand{\lem}{\begin{inspr}\begin{lemma}}
\newcommand{\elem}{\end{lemma}\end{inspr}}

\newcommand{\opm}{\begin{inspr}\begin{opmerking}}
\newcommand{\eopm}{\end{opmerking}\end{inspr}}

\newcommand{\voorw}{\begin{inspr}\begin{voorwaarde}}
\newcommand{\evoorw}{\end{voorwaarde}\end{inspr}}

\newcommand{\bew}{\vspace{-0.3cm}\begin{itemize}\item[ ] \bf Proof\rm: }
\newcommand{\ebew}{\hfill $\qed$ \end{itemize}}

\newcommand{\ssnl}{\vskip 4pt} 

\newcommand{\snl}{\vskip 7pt} 
\newcommand{\nl}{\vskip 12pt} 

\newcommand{\Cal}{\mathcal}

\newcommand{\rood}{\color{red}}
\newcommand{\blauw}{\color{blue}}
\newcommand{\magenta}{\color{magenta}}

\newcommand{\ot}{\otimes}
\newcommand{\di}{\diamond}

\newcommand{\tl}{\triangleleft}
\newcommand{\tr}{\triangleright}

\newcommand{\tussenen}{\qquad\quad\text{and}\qquad\quad}
\newcommand{\tussen}{\qquad\quad\qquad\quad}

\numberwithin{thm}{section}   
\numberwithin{equation}{section} 

\definecolor{bole}{rgb}{0.47, 0.27, 0.23}
\definecolor{auburn}{rgb}{0.43, 0.21, 0.1}
\definecolor{brown}{rgb}{0.59, 0.29, 0.0}

\begin{document}

\centerline{\bf \Large Weak Multiplier Hopf Algebras III}
\snl\centerline{\bf Integrals and duality}
\vspace{13pt}
\centerline{\it  Alfons Van Daele \rm $^{(1)}$ and \it Shuanhong Wang \rm $^{(2)}$}
\bigskip\bigskip
{\bf Abstract} 
\nl 
Let $(A,\Delta)$ be a {\it weak multiplier Hopf algebra} as introduced in \cite{VD-W4} (see also \cite{VD-W3}). It is a pair of a non-degenerate algebra $A$, with or without identity, and a coproduct $\Delta$ on $A$, satisfying certain properties. The main difference with multiplier Hopf algebras is that now, the canonical maps $T_1$ and $T_2$ on $A\ot A$, defined by
$$T_1(a\ot b)=\Delta(a)(1\ot b)
\qquad\quad\text{and}\qquad\quad
T_2(c\ot a)=(c\ot 1)\Delta(a),$$
are {\it no longer} assumed to be {\it bijective}. Also recall that a weak multiplier Hopf algebra is called {\it regular} if its antipode is a bijective map from $A$ to itself.
\snl
In this paper, we introduce and study the notion of {\it integrals} on such regular weak multiplier Hopf algebras. A left integral is a non-zero linear functional {\it on} $A$ that is left invariant (in an appropriate sense). Similarly for a right integral. 
\snl
For a  regular weak multiplier Hopf algebra $(A,\Delta)$ with (sufficiently many) integrals, we construct the dual $(\widehat A,\widehat\Delta)$. It is again a regular weak multiplier Hopf algebra with (sufficiently many) integrals. This duality extends the known duality of finite-dimensional weak Hopf algebras to this more general case. It also extends the duality of multiplier Hopf algebras with integrals, the so-called algebraic quantum groups. For this reason, we will sometimes call a regular weak multiplier Hopf algebra with enough integrals an {\it algebraic quantum groupoid}.
\snl
We discuss the relation of our work with the work on duality for algebraic quantum groupoids by Timmermann \cite{T2}.
\snl
We also illustrate this duality with a particular example in a separate paper (see \cite{VD8}). In this paper, we only mention the main definitions and results for this example.  
However, we do consider the two natural weak multiplier Hopf algebras associated with a groupoid in detail and show that they are dual to each other in the sense of the above duality.
\nl 
Date: {\it 25 September 2017}

\vspace{40pt}
\hrule
\smallskip
\begin{itemize}
\item[$^{(1)}$] Department of Mathematics, University of Leuven, Celestijnenlaan 200B, \newline
B-3001 Heverlee, Belgium. E-mail: \texttt{alfons.vandaele@kuleuven.be}
\item[$^{(2)}$] Department of Mathematics, Southeast University, Nanjing 210096, China. \newline 
E-mail: \texttt{shuanhwang@seu.edu.cn}
\end{itemize}
\newpage

\setcounter{section}{-1}  

\section{\hspace{-17pt}. Introduction} \label{s:introduction}  

With this paper, we continue our work on weak multiplier Hopf algebras as initiated in \cite{VD-W3}, \cite{VD-W4} and \cite{VD-W5}. 
\snl
Let $(A,\Delta)$ be {\it a regular weak multiplier Hopf algebra} as in Definition 4.1 of \cite{VD-W4}. The underlying algebra $A$ is a non-degenerate algebra over $\Bbb C$ and the coproduct is a homomorphism from $A$ to $M(A\ot A)$. It is {\it not} assumed to be non-degenerate, but there are certain constraints on the ranges and the kernels of the canonical maps $T_1$ and $T_2$ defined on $A\ot A$ as
\begin{equation*}
T_1(a\ot b)=\Delta(a)(1\ot b)
\qquad\quad\text{and}\qquad\quad
T_2(c\ot a)=(c\ot 1)\Delta(a).
\end{equation*}
More precisely, the images and the kernels of these maps are all characterized in terms of a single idempotent $E\in M(A\ot A)$. This idempotent is unique if it exists and it is called the canonical idempotent of $(A,\Delta)$.

\snl 
By assumption, there is a {\it counit} $\varepsilon$. This is a linear map from $A$ to $\mathbb C$. In general it is not a homomorphism. The coproduct is assumed to be full, i.e.\ the legs are all of $A$. Consequently the counit is unique. There is {\it a unique antipode} $S$. It is an anti-algebra map and an anti-coalgebra map. Because we assume that $(A,\Delta)$ is regular, it is a bijective map from $A$ to itself. 
\snl
The {\it source} and {\it target maps} $\varepsilon_s:A\to M(A)$ and $\varepsilon_t:A\to M(A)$ are defined by
\begin{equation*}
\varepsilon_s(a)=\sum_{(a)}S(a_{(1)})a_{(2)}
\tussenen
\varepsilon_t(a)=\sum_{(a)}a_{(1)}S(a_{(2)}).
\end{equation*}
We are using the Sweedler notation here (see further in this introduction under the item {\it Conventions and notations}).
It has been carefully argued in \cite{VD-W5} that these maps are  well-defined with values in $M(A)$.
\snl
The images $\varepsilon_s(A)$ and $\varepsilon_t(A)$ are subalgebras of $M(A)$. They are called the {\it source} and {\it target algebras}. 
They are commuting non-degenerate subalgebras of $M(A)$. In the regular case, they embed in $M(A)$ in such a way that their multiplier algebras $M(\varepsilon_s(A))$ and $ M(\varepsilon_t(A))$ still embed in $M(A)$. These multiplier algebras are denoted by $A_s$ and $A_t$ resp. They are still commuting subalgebras of $M(A)$. 
\snl
The source and target algebras $\varepsilon_s(A)$ and $\varepsilon_t(A)$ can be identified resp.\ with the left and the right leg of $E$. In fact we have $E\in M(\varepsilon_s(A)\ot \varepsilon_t(A))$.
\snl
Remark that in earlier papers on the subject, we 
called $A_s$ and $A_t$ the source and target algebras. However, in the second version of our second paper on the subject  \cite{VD-W5}, we changed the terminology. 
\snl
For more details about the short review over regular multiplier Hopf algebra given above, we refer to earlier work, in particular \cite{VD-W4} and \cite{VD-W5}.
\nl
\bf The examples associated with a groupoid \rm
\nl
Consider a groupoid $G$. One can associate two regular weak multiplier Hopf algebras.

\snl
{\it First} there is the algebra $A$, defined as the space $K(G)$ of complex functions on $G$ with finite support and pointwise product. The coproduct $\Delta$ on $K(G)$ is defined by  
\begin{equation*}
\Delta(f)(p,q)=
\begin{cases}
		f(pq) & \text{if $pq$ is defined},\\
				0 & \text{otherwise}.
\end{cases}
\end{equation*}

\snl
In this example, the algebra $A_s$
 is the algebra of functions on $G$ so that $f(p)=f(q)$ whenever $p,q\in G$ satisfy $s(p)=s(q)$ (where $s$ is the source map of the groupoid). The algebra $A_t$ 
 consists of functions $f$ on $G$ so that $f(p)=f(q)$ if $t(p)=t(q)$ for $p,q\in G$ (where $t$ is the target map). 
\snl
For the {\it second case}, we take the algebra $B$, defined as the groupoid algebra $\Bbb C G$ of $G$. If we use $p\mapsto \lambda_p$ for the canonical embedding of $G$ in $\Bbb C G$, then for $p,q\in G$ we have $\lambda_p\lambda_q=\lambda_{pq}$  if $pq$ is defined and $\lambda_p\lambda_q=0$ otherwise. 
The coproduct is defined by $\Delta(\lambda_p)=\lambda_p\ot \lambda_p$ for all $p\in G$.
Here the multiplier algebras $B_s$ and $B_t$ of the source and target algebras coincide and it is the multiplier algebra of the span of elements of the form $\lambda_e$ where $e$ is a unit of $G$.
\snl

Remark that the first algebra has no identity if $G$ is infinite. The second one has no identity if the set of units of $G$ is infinite. In particular, only when there is an identity, they are weak Hopf algebras. If not, they are weak {\it multiplier} Hopf algebras. 
\snl
However, in the two cases local units always exist. This is true by a general result, but of course, it is easy to find them for these two examples. 
\snl
We refer to Examples 1.15 and 1.16 in \cite{VD-W4} and Example 3.1 in \cite{VD-W5} for more details on these examples. 
\snl
These two cases are {\it dual} to each other. The duality is given by $\langle f,\lambda_p\rangle=f(p)$ whenever $f\in K(G)$ and $p\in G$. We will give more details about this duality in Section \ref{s:examples} of this paper where we illustrate the duality for regular weak multiplier Hopf algebras with integrals as developed in Section \ref{s:duality}.
\nl
Let us now look at {\it integrals} for these two examples. 
\snl
First take for $A$ the algebra $K(G)$. Let $g$ be any complex function on $G$ satisfying $g(q)=g(q')$ whenever $q$ and $q'$ have the same source (i.e.\ we have $g\in A_s$). Define a linear functional $\varphi$ on $A$ by $\varphi(f)=\sum_q g(q)f(q)$. Then $\varphi$ will be a left integral and all left integrals are of this form. We refer to Section \ref{s:integrals} for the precise definition of a left integral and to Section \ref{s:examples} for the proof of this result. Similarly, a right integral $\psi$ is a linear functional on $A$ of the form $\psi(f)=\sum_q h(q)f(q)$ where now $h$ is a function such that $h(q)=h(q')$ whenever $q$ and $q'$ have the same target, i.e.\ when $h\in A_t$. Again we refer to Section \ref{s:examples} for a proof of this result. It can happen that $\varphi$ will be left and right invariant (e.g.\ if $g=h=1$) but not all left integrals will also be right integrals.
\snl
Next consider the groupoid algebra $\Bbb C G$. Let $g$ be a complex function on $G$ with support in the set of units. Define $\varphi(\lambda_p)=g(p)$ for all $p$. This will give a left integral and all left integrals are of this form. Now the left and right integrals are the same because the coproduct is coabelian.
\snl
Again, we refer to Section \ref{s:examples} where we will treat these examples in detail. 
\nl
What we study in this paper generalizes the results on integrals, described above for the regular weak multiplier Hopf algebras $K(G)$ and $\Bbb C G$, associated with a groupoid $G$, to general regular weak multiplier Hopf algebras. 
\nl
\bf Relation with the work of Timmermann on duality of algebraic quantum groupoids \rm \cite{T2} 
\nl
While preparing this manuscript, we came across recent work by  Timmermann. He studies integrals and duality for regular multiplier Hopf algebroids in \cite{T1} and \cite{T2}. As an application, based on the notion of integrals on weak multiplier Hopf algebras as defined in \cite{K-VD}, together with the characterization of regular multiplier Hopf algebroids in terms of weak multiplier Hopf algebras, studied in \cite{T-VD1}, he obtains also a duality theorem for regular weak multiplier Hopf algebras, similar as the one obtained in this paper.
\snl
The relation of his work and this one will be discussed further, and in detail, at several places in this paper. 
To summarize, we have added a couple of remarks in Section \ref{s:conclusions}. 
\snl
We feel it is appropriate to mention that the present paper has a long history. As mentioned already, it is a natural continuation of our previous work on the subject (\cite{VD-W3}, \cite{VD-W4} and \cite{VD-W5}). 
Moreover, the main results included in this paper have been presented at several occasions in the past (see e.g. \cite{VD9} and \cite{VD10}, and also \cite{VD-Wo}). Finally, integrals on weak multiplier Hopf algebras are also studied in \cite{K-VD}.
\nl
\bf Content of the paper.\rm
\nl
In {\it Section} \ref{s:integrals}, we introduce the notion of left and right integrals on regular weak multiplier Hopf algebras. We discuss the definition and compare with the literature on weak Hopf algebras and other cases where integrals are studied. 
\snl 
To do this, we first prove some basic properties  and use them to give alternative definitions of left and right invariance of linear functionals. Then we prove properties that are very similar as in the case of multiplier Hopf algebras with integrals. However, we state  clearly what the main differences are. We do not have automatically that a non-zero integral is faithful and so in order to be able to get nice results, we need to assume that there are {\it enough} integrals. If there is a single faithful integral, we get the scaling constant, the modular element as well as the modular automorphisms, just as for multiplier Hopf algebras. 
In the general case, one can probably get partial results along these lines. We refer to Section \ref{s:conclusions} where we discuss this further.
\snl

In this section we already define the dual space $\widehat A$ and give some characterizations of it. However, it is only in Section \ref{s:duality} that we make this dual space again into a weak multiplier Hopf algebra. The basic examples are considered in Section \ref{s:examples}. 
\snl
In {\it Section} \ref{s:duality} we study duality. We show that the dual $\widehat A$ of $A$, already introduced as a subspace of the dual space $A'$ in Section \ref{s:integrals}, can be equipped with a product and  a coproduct, adjoint to the coproduct and product on $A$ respectively. This turns $\widehat A$ again into a regular weak multiplier Hopf algebra.
\snl
The construction of the coproduct on the dual is not so obvious. We show that the adjoints of the various canonical maps and their generalized inverses of $A$ can be found on the dual $\widehat A$. This provides the coproduct on $\widehat A$. We show that indeed, this is adjoint to the product on $A$ by first extending the pairings to the multiplier algebras (on one side).
\snl
We obtain explicit formulas for the dual objects: the counit, the canonical idempotent, the antipode and the source and target maps and source and target algebras. In particular, we find that the source algebra of $\widehat A$ is isomorphic with the target algebra of $A$ and that the target algebra of $\widehat A$ is isomorphic with the source algebra of $A$. This is known for weak Hopf algebras, see e.g.\ Lemma 2.6 in \cite{B-N-S}.
\snl
Finally, we construct integrals on the dual and show that there is a faithful set of integrals, also on the dual. 
If the orginal has a single faithful integral, the same is true for the dual.
\snl
It follows that the dual of $(\widehat A,\widehat \Delta)$ can be constructed. It is canonically isomorphic with the original weak multiplier Hopf algebra $(A,\Delta)$. 
\nl
In {\it Section} \ref{s:examples} we discuss 
examples.
\snl
First we look at the two weak multiplier Hopf algebras associated with a groupoid $G$. The first one is the function algebra $K(G)$ while the second one is the groupoid algebra $\mathbb C G$. The two weak multiplier Hopf algebras have been discussed in earlier papers, see e.g.\ Examples 1.15 and 1.16 in \cite{VD-W4}. Here we obtain the integrals. It turns out that in the two cases, there is a faithful left integral that is also a right integral. We also show in detail how the groupoid algebra is the dual of the function algebra, hereby applying the general duality theory as developed in Section \ref{s:duality}.
\snl
We also briefly treat the weak multiplier Hopf algebra associated with a separability idempotent. This example is already considered in Proposition 2.8 of \cite{VD-W5}. Here, we find the integrals on this weak multiplier Hopf algebra and again, we follow the results of Section \ref{s:duality} to obtain the dual of this weak multiplier Hopf algebra. We only formulate the main results about this case. We refer to a separate paper (\cite{VD8}) where we treat this example in great detail. By doing so, we illustrate various aspects of the theory as developed in this paper. 
Remark that this example is also treated in \cite{T2}
\nl
In {\it Section} \ref{s:conclusions} we draw some conclusions and discuss possible further research. 
\nl
\bf Conventions and notations \rm
\snl
We only work with algebras $A$ over $\Bbb C$ (although we believe that this is not essential). We do not assume that they are unital but we need that the product is non-degenerate. Our algebras are all idempotent (that is $A^2=A$). In some situations, this is a condition, while at other places, it follows from the other axioms. In certain parts of the paper, we even need the algebras to have local units. Also this is sometimes a condition, while in other cases, it is a consequence. Then of course, the product is automatically non-degenerate and also the algebra is idempotent.
\snl
When $A$ is such an algebra, we use $M(A)$ for the multiplier algebra of $A$. When $m$ is in $M(A)$, then by definition we  have elements  $am$ and $mb$ in $A$ for all $a,b\in A$ and $(am)b=a(mb)$. The algebra $A$ sits in $M(A)$ as an essential two-sided ideal and $M(A)$ is the largest algebra with identity having this property. 
\snl
We consider $A\ot A$, the tensor product of $A$ with itself. It is again an idempotent, non-degenerate algebra and we can consider the multiplier algebra $M(A\ot A)$. The same is true for a multiple tensor product. We will no longer use $\sigma$ for the flip map on $A\ot A$ (as we have done sometimes in earlier papers) because this will be used now for the modular automorphism. In stead, we will use $\zeta$ for the flip map on $A\ot A$, as well as for its natural extension to $M(A\ot A)$.
\snl
We use $1$ for the identity in any of these  multiplier algebras. On the other hand, we mostly use $\iota$ for the identity map on $A$ (or other spaces), although sometimes, we also write $1$ for  this map. The identity element in a group is denoted by $e$. If $G$ is a groupoid, we will also use $e$ for units. Units are considered as being elements of the groupoid and we use $s$ and $t$ for the source and target maps from $G$ to the set of units. 
\snl
The space of all linear functionals on $A$ is denoted by $A'$. A linear functional $\omega$ on an algebra $A$ is called {\it faithful} if given $a\in A$ we have $a=0$ if either $\omega(ab)=0$ for all $b\in A$ or $\omega(ba)=0$ for all $b\in A$. 
\snl
When $A$ is an algebra, we denote by $A^{\text{op}}$ the algebra obtained from $A$ by reversing the product. When $\Delta$ is a coproduct on $A$, we denote by $\Delta^{\text{cop}}$ the coproduct on $A$ obtained by composing $\Delta$ with the flip map $\zeta$.
\snl
For a coproduct $\Delta$, as we define it in Definition 1.1 of \cite{VD-W4}, we assume that $\Delta(a)(1\ot b)$ and $(c\ot 1)\Delta(a)$ are in $A\ot A$ for all $a,b,c\in A$. This allows us to make use of the {\it Sweedler notation} for the coproduct. The reader who wants to have a deeper understanding of this, is referred to \cite{VD7} where the use of the Sweedler notation for coproducts that do not map into the tensor product, but rather in its multiplier algebra, is explained in detail. 
\nl
\bf Basic references \rm
\nl
For the theory of Hopf algebras, we refer to the standard works of Abe \cite{A} and Sweedler\cite{S}. For multiplier Hopf algebras and integrals on multiplier Hopf algebras, we refer to \cite{VD1} and \cite{VD2}. Weak  Hopf algebras have been studied in \cite{B-N-S} and \cite{B-S} and more results are found in \cite{Ni} and \cite{N-V1}. Various other references on the subject can be found in \cite{Va}. In particular, we refer to \cite{N-V2} because we will use notations and conventions from this paper when dealing with weak multiplier Hopf algebras.
\snl
Weak multiplier Hopf algebras have been introduced in \cite{VD-W4}. See also \cite{VD-W3} for a preliminary paper on the subject. The source and target algebras and the source and target maps for a weak multiplier Hopf algebra are studied in detail in \cite{VD-W5}. In the more recent second version of the paper, more results are obtained for possibly non-regular weak multiplier Hopf algebras.  Integrals on regular weak multiplier Hopf algebras are also considered in \cite{K-VD} where a form of the Larson-Sweedler theorem is obtained in the framework of regular weak multiplier Hopf algebras. 
\snl
Finally, for the theory of groupoids, we refer to \cite{Br}, \cite{H},  \cite{P} and  \cite{R}. 
\nl
\bf Acknowledgements \rm
\nl
The first named author (Alfons Van Daele) would like to express his thanks to his coauthor Shuanhong Wang for motivating him to start the research on weak multiplier Hopf algebras and for the
hospitality when visiting the University of Nanjing in 2008 for the first time and later again in
2012, 2014 and 2016.
\snl
The second named author (Shuanhong Wang) would like to thank his coauthor for his help and
advice when visiting the Department of Mathematics of the University of Leuven in Belgium during
several times the past years. 
\snl
This work is partially supported by the NSF of China (No 11371088)  and the NSF of Jiangsu Province (No. BK20171348).
\snl
We would also like to thank Thomas Timmermann for discussions about the relation of this work and his recent work on duality for algebraic quantum groupoids \cite{T2}. 
\nl\nl

\section{\hspace{-17pt}. Integrals on weak multiplier Hopf algebras}\label{s:integrals} 

In this section, we treat integrals on regular weak multiplier Hopf algebras. We will make a remark about the theory of integrals on possibly non-regular weak multiplier Hopf algebras in Section \ref{s:conclusions}.
\snl
Integrals on weak multiplier bialgebras have been studied before in the context of a version of the Larson-Sweedler theorem for weak multiplier Hopf algebras, see \cite{K-VD}. In the context of multiplier Hopf algebroids, integrals were defined and investigated in \cite{T1}. Here, we give a self-contained treatment. 
\snl
We fix a {\it regular} weak multiplier Hopf algebra $(A,\Delta)$. We start with the definition of a left, respectively a right integral on $A$. We first prove some basic properties that give equivalent definitions for left and right invariance of linear functionals. We compare our definition with the ones found in the literature for weak Hopf algebras (and other cases).
\snl
Then we develop the theory of integrals further. We obtain results that are very similar as those for multiplier Hopf algebras. The main difference is that integrals are not necessarily unique and they are not automatically faithful. Therefore, to obtain a useful theory, we will need to assume not only the existence of integrals, but we have to require that there are enough integrals, in a sense to be explained.
\nl
\bf Left and right integrals \rm
\nl
Before we give the main definitions, recall that for any element $a\in A$ and any linear functional $\omega$ on $A$, we can define a multiplier $x\in M(A)$ by
\begin{equation*}
	xb=(\iota\ot\omega)(\Delta(a)(b\ot 1)) 
	\tussenen
	bx=(\iota\ot\omega)(b\ot 1)\Delta(a))
\end{equation*}
where $b$ is in $A$. This multiplier $x$ is denoted as $(\iota\ot\omega)\Delta(a)$. Similarly we can define $(\omega\ot\iota)\Delta(a)$ in $M(A)$ for $a\in A$ and any linear functional $\omega$ on $A$. We need regularity of the coproduct to do this but this is satisfied as we work with a regular weak multiplier Hopf algebra.
\snl
The above property is used in the following definition.

\defin\label{defin:1.1}
A linear functional $\varphi:A\to \Bbb C$ is called {\it left invariant} if $(\iota\ot\varphi)\Delta(a)\in A_t$ for all $a\in A$. Similarly, a linear functional $\psi$ on $A$ is called {\it right invariant} if $(\psi\ot\iota)\Delta(a)\in A_s$ for all $a\in A$. A non-zero left invariant functional is called a {\it left integral} and a non-zero right invariant functional is called a {\it right integral} on the weak multiplier Hopf algebra $(A,\Delta)$.
\edefin

Recall that $A_t$ and $A_s$ are  subspaces of $M(A)$ and so the above definition makes sense.

\opm
i) If $(A,\Delta)$ is actually a multiplier Hopf algebra, so that $A_t$ and $A_s$ are nothing else but scalar multiples of the identity, it follows from the requirement in the definition, that for any $a\in A$ we have $(\iota\ot\varphi)\Delta(a)=\varphi(a)1$ when $\varphi$ is left invariant and $(\psi\ot\iota)\Delta(a)=\psi(a)1$ when $\psi$ is right invariant. So, the above definition is consistent with the definition of integrals in the case of a regular multiplier Hopf algebra.
\snl
ii) As the antipode $S$ flips the coproduct and maps $A_t$ to $A_s$ and vice versa, we will have that $\varphi\circ S$ is a right integral when $\varphi$ is a left integral and that $\psi\circ S$ is a left integral when $\psi$ is a right integral. Also $\varphi\circ S^2$  will be a left integral when $\varphi$ is a left integral. However, as we no longer have uniqueness of integrals, we can not expect that this is a scalar multiple of $\varphi$. We will see later how this property is modified as a consequence of Proposition \ref{prop:2.7}, in the case of a faithful integral.
\eopm

If $(A,\Delta)$ is a weak Hopf algebra, then we have $\varepsilon_t(x)=x$ for all $x\in A_t$ and so for all $a$ we have $(\iota\ot\varphi)\Delta(a)=(\varepsilon_t\ot\varphi)\Delta(a)$ when $\varphi$ is left invariant. This shows that our definition of a left integral coincides with that of a left integral on a weak Hopf algebra as it is found in \cite{B-N-S}. Similarly for a right integral. 
\snl
However, we cannot simply import the definition of integrals as defined for weak Hopf algebras in this setting. It would not make sense to write $(\iota\ot\varphi)\Delta(a)=(\varepsilon_t\ot\varphi)\Delta(a)$ because we have not defined $\varepsilon_t$ on the multiplier algebra. Nevertheless, this just turns out to be a matter of interpretation as we will see in the next proposition.

\prop\label{prop:2.3}
Let $\varphi$ be a left integral on $A$. Then 
\begin{equation}
 	(\iota\ot\varphi)\Delta(a)= \sum_{(a)} a_{(1)} S(a_{(2)})\varphi(a_{(3)}) \label{eqn:2.1}
\end{equation}
for all $a\in A$. Similarly, when $\psi$ is a right integral, we have
\begin{equation}
	(\psi\ot\iota)\Delta(a)=\sum_{(a)} \psi(a_{(1)}) S(a_{(2)}) a_{(3)}\label{eqn:2.2}
\end{equation}
for all $a$.
\eprop

\bew
The formula (\ref{eqn:2.1}) is given a meaning if we multiply with an element $c$ of $A$ from the left. The left hand side becomes an element in $A$ and this is true also for the right hand side. Indeed, the element $c$ will first cover $a_{(1)}$ and subsequently, also $a_{(2)}$ will be covered. See e.g.\ Remark 1.2 in \cite{VD-W5} about coverings in these formulas.  The formula (\ref{eqn:2.2}) is given a meaning if we multiply with an element of $A$ from the right.
\ssnl
To prove the formula (\ref{eqn:2.1}) take $a,b,c\in A$ and let $x=(\iota\ot\varphi)\Delta(a)$. By definition, we know that $x\in A_t$ and therefore, by a property of the target map (see Proposition 2.7 in \cite{VD-W5}), we have that $x\varepsilon_t(b)=\varepsilon_t(xb)$. We multiply this equality with $c$ from the left and use the definition of $\varepsilon_t$. This gives
\begin{align*}
	cx\varepsilon_t(b)=c\varepsilon_t(xb)
	&=\sum_{(a),(b)}ca_{(1)}b_{(1)}S(a_{(2)}b_{(2)})\varphi(a_{(3)})\\
	&=\sum_{(a),(b)}ca_{(1)}b_{(1)}S(b_{(2)})S(a_{(2)})\varphi(a_{(3)})\\
	&=\sum_{(a)}ca_{(1)}\varepsilon_t(b)S(a_{(2)})\varphi(a_{(3)}).
\end{align*}
Now think of $ca_{(1)}\ot a_{(2)}$ as $p\ot q$ in $A\ot A$.
For the right hand side of the last expression above we then consider
\begin{equation*}
	p\varepsilon_t(b)S((\iota\ot\varphi)\Delta(q)).
\end{equation*} 
As $\varphi$ is left invariant, we know that $(\iota\ot\varphi)\Delta(q)$ is in $A_t$ and so $S((\iota\ot\varphi)\Delta(q))$ will be in $A_s$. Now we use that $A_t$ and $A_s$ commute. This implies that 
\begin{equation*}
	p\varepsilon_t(b)S((\iota\ot\varphi)\Delta(q))=pS((\iota\ot\varphi)\Delta(q))\varepsilon_t(b).
\end{equation*}
Replace again $p\ot q$ by $ca_{(1)}\ot a_{(2)}$ and insert the result in the original formula above. This will give
\begin{equation*}
	cx\varepsilon_t(b)=\sum_{(a)}ca_{(1)}S(a_{(2)})\varepsilon_t(b)\varphi(a_{(3)}).
\end{equation*}
As this is true for all $b$, we find 
$cx=\sum_{(a)}ca_{(1)}S(a_{(2)})\varphi(a_{(3)})$
for all $c$ (as a consequence of Proposition 2.10 in \cite{VD-W5}). This gives the correct interpretation of the formula (\ref{eqn:2.1}) as we want to prove it.
\ssnl
The argument for the right integral is completely similar.
\ebew

The formulas (\ref{eqn:2.1}) and (\ref{eqn:2.2}) in the previous proposition suggest that $(\iota\ot\varphi)\Delta(a)\in \varepsilon_t(A)$ when $\varphi$ is a left integral and that $(\psi\ot\iota)\Delta(a)\in \varepsilon_s(A)$ when $\psi$ is a right integral. We will see shortly that this is indeed the case. However, we can not conclude it from the above formulas.
\snl
These results are intimately related with the following result obtained in \cite{K-VD}.  Recall that $E$ is the canonical idempotent of $(A,\Delta)$.

\prop\label{prop:2.4a}
Denote 
\begin{align*}
&F_1=(\iota\ot S)E \tussenen F_3=(\iota\ot S^{-1})E \\
&F_2=(S\ot \iota)E \tussenen F_4=(S^{-1}\ot \iota)E.
\end{align*}	
Then, if $\varphi$ is a left integral and if $\psi$ is a right integral, we have for all $a$ in $A$,
\begin{align}
(\iota\ot\varphi)\Delta(a) &=(\iota\ot\varphi)(F_2(1\ot a))=(\iota\ot\varphi)((1\ot a)F_4) \label{eqn:2.3}\\
(\psi\ot\iota)\Delta(a)&=(\psi\ot\iota)((a\ot 1)F_1)=(\psi\ot\iota)(F_3(a\ot 1)).\label{eqn:2.4}
\end{align}
\eprop

The results are found in Proposition 2.7 of \cite{K-VD} and the proofs can be found in that paper. The context of \cite{K-VD} is a little different, but this has no effect for the proofs.
\snl
In fact, the formulas are reformulations of the result in Proposition \ref{prop:2.3}. 
\snl
Consider e.g.\ the formula
\begin{equation*}
E(1\ot a)=\sum_{(a)}a_{(2)}S^{-1}(a_{(1)})\ot a_{(3)}
\end{equation*}
for $a\in A$ as found from formula (1.5) in Section 1 of \cite{VD-W5}. Apply $S$ on the first factor to get
\begin{equation*}
F_2(1\ot a)=\sum_{(a)}a_{(1)}S(a_{(2)})\ot a_{(3)}
\end{equation*}
and apply $\varphi$ on the second factor. We get
\begin{equation*}
(\iota\ot\varphi)(F_2(1\ot a))=\sum_{(a)}a_{(1)}S(a_{(2)})\varphi(a_{(3)}).
\end{equation*}
We see that this is the formula (\ref{eqn:2.1}) in Proposition \ref{prop:2.3} above. Similarly we have 
\begin{equation*}
(a\ot 1)E=\sum_{(a)}a_{(1)}\ot S^{-1}(a_{(3)})a_{(2)}
\end{equation*}
as in formula (1.6) in Section 1 of \cite{VD-W5}. Now apply $S$ on the second factor and $\psi$ on the first one to arrive at
\begin{equation*}
(\psi\ot \iota)((a\ot 1)F_1)=\sum_{(a)}\psi(a_{(1)}) S(a_{(2)})a_{(3)}
\end{equation*}
giving formula (\ref{eqn:2.2}) of Proposition \ref{prop:2.3} above.
\snl
If we write down these formulas for $(A^{\text{op}},\Delta)$ we find the two other formulas. Indeed, by replacing the product in $A$ by the opposite product, there is no effect on the coproduct, so that $E$ remains the same, as well as left and right invariance. Only the antipode has to be replaced by it inverse. So, $F_1$ and $F_2$ become $F_3$ and $F_4$ and taken into account that we reverse the order, we indeed get the other two formulas.
\snl
Remark that in Section 5 of \cite{T2}, the first equality in Equation (\ref{eqn:2.3}) is used as a definition for left invariance while the first equality in Equation (\ref{eqn:2.4}) as a definition for right invariance. In this paper, we have another definition of left and right invariance (see Definition \ref{defin:1.1}) while the equations used in Definition 5.1.3 of \cite{T2} are properties. 

\nl
We will now continue the study of left and right integrals on a regular weak multiplier Hopf algebra. We aim at similar properties as for regular multiplier Hopf algebras. The following result is, surprisingly enough, precisely the same as in that case.

\prop\label{prop:2.4}
Let $\varphi$ be a linear functional on $A$. Given $a,b\in A$, define
\begin{equation}
	c=(\iota\ot\varphi)(\Delta(a)(1\ot b))
	\tussenen
	d=(\iota\ot\varphi)((1\ot a)\Delta(b)).\label{eqn:2.5}
\end{equation}
These elements belong to $A$. Then $\varphi$ is left invariant if and only if we have $d=S(c)$ for all $a$ and $b$. Similarly, let $\psi$ be a linear functional on $A$ and now define, given $a,b\in A$ the elements $c,d\in A$ by
\begin{equation}
	c=(\psi\ot\iota)((a\ot 1)\Delta(b))
	\tussenen
	d=(\psi\ot\iota)(\Delta(a)(b\ot 1)).\label{eqn:2.6}
\end{equation}
Then $\psi$ is right invariant if and only if $d=S(c)$ for all $a$ and $b$.
\eprop

\bew
i) First assume that $\varphi$ is left invariant. Take any $a$ and $b$ in $A$ and define the elements $c$ and $d$ as in the formulation of the proposition. We have
\begin{equation}
	(1\ot a)\Delta(b)=\sum_{(a)} (S(a_{(1)}) \ot 1)\Delta(a_{(2)}b)\label{eqn:1.7}
\end{equation}
and if we apply $\varphi$ on the second factor, we find
\begin{equation*}
	d=(\iota\ot\varphi)((1\ot a)\Delta(b))=\sum_{(a)} S(a_{(1)}) ((\iota\ot\varphi)\Delta(a_{(2)}b)).
\end{equation*}
Now apply $\Delta$. Because  $\Delta(px)=\Delta(p)(x\ot 1)$ when $p\in A$ and $x\in A_t$, it follows from the left invariance of $\varphi$ that
\begin{equation*}
	\Delta(d)=\sum_{(a)} \Delta(S(a_{(1)}))((\iota\ot\varphi)\Delta(a_{(2)}b)\ot 1).
\end{equation*}
Next multiply with an element $q$ from the right in the second factor. If we use that $S$ flips the coproduct, we find
\begin{align*}
	\Delta(d)(1\ot q) 
	&=\sum_{(a)} (S(a_{(2)})\ot S(a_{(1)})q)((\iota\ot\varphi)\Delta(a_{(3)}b)\ot 1) \\
	&=\sum_{(a)} S(a_{(2)})((\iota\ot\varphi)\Delta(a_{(3)}b)) \ot S(a_{(1)})q.
\end{align*}
We use the formula  (\ref{eqn:1.7}) with $a$ replaced by $a_{(2)}$ and we find
\begin{equation*}
	\Delta(d)(1\ot q)=\sum_{(a)} (\iota\ot\varphi)((1\ot a_{(2)})\Delta(b)) \ot S(a_{(1)})q.
\end{equation*}
We can now safely apply $\varepsilon$ on the first leg of this equality and we get precisely $dq=S(c)q$ where $c=\sum_{(a)} a_{(1)} \varphi(a_{(2)}b)$. We can cancel $q$ and this proves one part of the proposition.
\ssnl
ii) Conversely, assume that $\varphi$ is a linear functional satisfying the requirement in the proposition. Then for all $p,q$ in $A$ we have
\begin{align*}
	(\iota\ot\varphi)\Delta(pq)
	&=\sum_{(p)} p_{(1)} (\iota\ot\varphi)((1\ot p_{(2)})\Delta(q))\\
	&=\sum_{(p)} p_{(1)} S((\iota\ot\varphi)(\Delta(p_{(2)})(1\ot q)))\\
	&=\sum_{(p)} p_{(1)} S(p_{(2)})\varphi(p_{(3)}q).
\end{align*}
We can cover these expressions if we multiply with an element of $A$  from the left. 
As $A^2=A$, this proves that $(\iota\ot\varphi)\Delta(a)\in \varepsilon_t(A)$ for all $a\in A$. And because $\varepsilon_t(A)$ is a subset of $A_t$, we have shown that $\varphi$ is left invariant.
\ssnl
iii) The proof for right invariant functionals is completely similar.
\ebew

We have made use of the Sweedler notation but most of the time, we have used proper coverings. At other places where covering is necessary, we have indicated how the expressions can be covered. The reader can verify that this can easily be done.
\snl
In the case of a weak Hopf algebra, that is when $A$ has an identity, we can take $a=1$ in Equation (\ref{eqn:2.5}) and we will obtain the first equality of Equation (\ref{eqn:2.3}). Similarly, if we take $b=1$ in Equation (\ref{eqn:2.6}) we will get the first equation in (\ref{eqn:2.4}). Doing this with $b=1$ in the first case and with $a=1$ in the second case, we get the other equalities.
\snl
At the end of the proof of the previous proposition, we have found that 
\begin{equation*}
	(\iota\ot\varphi)(\Delta(pq))=\sum_{(p)} p_{(1)} S(p_{(2)})\varphi(p_{(3)}q)=\sum_{(p)} \varepsilon_t(p_{(1)})\varphi(p_{(2)}q)
\end{equation*}
for all $p\in A$ when $\varphi$ is a left integral. If the algebra has an identity, as in the case of weak Hopf algebras, we can take $q=1$ and this will give again the formula obtained in Proposition \ref{prop:2.3}. In general, it follows that $(\iota\ot\varphi)(\Delta(a))\in \varepsilon_t(A)$ when $a\in A$. This result was expected but not yet proven (see a remark following Proposition \ref{prop:2.3}). Of course, when $\psi$ is a right integral, we will have $(\psi\ot\iota)\Delta(a)\in \varepsilon_s(A)$.
\snl
Remark in passing that in the theory of {\it algebraic quantum hypergroups} (see e.g.\ \cite{D-VD}), the above formulas are still true for left and right integrals. In fact, they are used to define integrals in that theory. In the theory of algebraic quantum hypergroups, it is not assumed that the coproduct is an algebra homomorphism.
\snl
Having exactly the same characterizing formulas for invariant functionals (Proposition \ref{prop:2.4}), one can expect that many of the other properties of integrals, as we know them in the theory of multiplier Hopf algebras and algebraic quantum hypergroups, will also be true here. We will see that this is the case, except for the uniqueness and the related property of faithfulness.
\snl
The claim above is nicely illustrated in the following proposition where we prove a set of formulas that will be crucial for the various properties of the integrals, as well as for the definition of the dual $\widehat A$ and to show that it is again a regular weak multiplier Hopf algebra (cf.\ Section \ref{s:duality}).

\prop\label{prop:2.5}
Let $\varphi$ be any left integral and $\psi$ any right integral on $A$. In the following formulas, $p$ and $q$ are arbitrary elements in $A$.
\snl
i) We have $\psi(xa)=\varphi(xb)$ for all $x\in A$ if
\begin{align*}
	a&= (\iota\ot\varphi)((\iota\ot S)(\Delta(p))(1\ot q))\\
	b&= (\psi\ot\iota)((S^{-1}\ot\iota)(\Delta(q))(p\ot 1)).
\end{align*}
ii) We have $\psi(ax)=\varphi(bx)$ for all $x\in A$ if
\begin{align*}
	a&=(\iota\ot\varphi)((1\ot q)(\iota\ot S^{-1})\Delta(p))\\
	b&= (\psi\ot\iota)((p\ot 1)(S\ot\iota)\Delta(q)).
\end{align*}
iii) We have $\psi(xa)=\varphi(bx)$ for all $x\in A$ if
\begin{align*}
	a&= (\iota\ot\varphi)((1\ot q)(\iota\ot S)\Delta(p))\\
	b&= (\psi\ot\iota)((S\ot\iota)(\Delta(q))(p\ot 1)).
\end{align*}
iv) We have $\psi(ax)=\varphi(xb)$ for all $x\in A$ if
\begin{align*}
	a&= (\iota\ot\varphi)((\iota\ot S^{-1})(\Delta(p))(1\ot q)) \\
	b&= (\psi\ot\iota)((p\ot 1)(S^{-1}\ot\iota)\Delta(q)).
\end{align*}
\eprop

\bew
To prove i), take $p,q\in A$ and also $x\in A$. Consider the expression \newline
$(\psi\ot\varphi)((\Delta(x)(p\ot q))$. 
Using the formula obtained in Proposition \ref{prop:2.4} for $\varphi$ with $a,b$ replaced by $p,q$, we find that 
\begin{align*}
	(\psi\ot\varphi)(\Delta(x)(p\ot q))
	&=\psi(\,\cdot\,p)((\iota\ot\varphi)(\Delta(x)(1\ot q))) \\
	&=\psi(\,\cdot\,p)(S^{-1}((\iota\ot\varphi)((1\ot x)\Delta(q)))) \\
	&=(\psi\ot\varphi)((1\ot x)(S^{-1}\ot\iota)(\Delta(q))(p\ot 1))\\
	&=\varphi(xb)
\end{align*}
where $b= (\psi\ot\iota)((S^{-1}\ot\iota)(\Delta(q))(p\ot 1))$. On the other hand, if we use the formula in Proposition \ref{prop:2.4} for $\psi$, we find that
\begin{align*} 
	(\psi\ot\varphi)(\Delta(x)(p\ot q))
	&=\varphi(\,\cdot\,q)((\psi\ot\iota)(\Delta(x))(p\ot 1)) \\
	&=\varphi(\,\cdot\,q)(S((\psi\ot\iota)(x\ot 1)\Delta(p))) \\
	&=(\psi\ot\varphi)((x\ot 1)(\iota\ot S)(\Delta(p))(1\ot q))\\
	&=\psi(xa)
\end{align*}
where $a= (\iota\ot\varphi)((\iota\ot S)(\Delta(p))(1\ot q))$. This proves i).
\snl
To prove ii) we start with the expression $(\psi\ot\varphi)((p\ot q)\Delta(x))$ and proceed as above. For the proof of iii), we take $(\psi\ot\varphi)((1\ot q)\Delta(x)(p\ot 1)$ and for iv) we finally start with $(\psi\ot\varphi)((p\ot 1)\Delta(x)(1\ot q))$. And also in these cases, we use the two formulas of Proposition \ref{prop:2.4}. 
\ebew

At the end of this section, we will compare our notions and results with the ones in the paper by Timmermann \cite{T2}. We will then also prove some more results on integrals as defined in this paper. However, these results are of minor importance for our approach.
\snl
First we consider the following special case.

\nl
\bf  The case of a single faithful integral \rm
\nl
In the event that there exists a single faithful integral, we find some nice consequences of these results. 
\snl
First observe the following. Assume that $\varphi$ is just a faithful linear functional. We claim that all elements in $A$ are linear combinations of elements of the form $a$ in each of the items in the previous proposition. 
\snl
Indeed, suppose e.g.\ that $\omega$ is a linear functional on $A$ that vanishes on all elements of the form
\begin{equation*}
 (\iota\ot\varphi)((\iota\ot S)(\Delta(p))(1\ot q))
\end{equation*}
as considered in item i) of the proposition. Replace $q$ by $qq'$ and write 
\begin{equation*}
x=(\omega\ot\iota)((\iota\ot S)(\Delta(p))(1\ot q)).
\end{equation*}
Then  $\varphi(xq')=0$ for all $q'$ and because $\varphi$ is assumed to be faithful, it follows that $x=0$. This is true for all $p,q$ and because the coproduct is full, we must have $\omega=0$. This proves the claim. 
Recall that fullness of the coproduct means that the legs are all of $A$.
\ssnl
A similar argument can be used for the three other cases.
\snl 
As a first consequence, we get the KMS property.

\prop\label{prop:2.7a}
If $\varphi$ is a faithful left integral, there exists an automorphism $\sigma$ of $A$ satisfying $\varphi(ab)=\varphi(b\sigma(a))$ for all $a,b$. It leaves $\varphi$ invariant.
\eprop

\bew
Because there is a faithful left integral $\varphi$, there is also a faithful right integral $\psi$. 
If we combine the formulas in the items i) and iii) of Proposition \ref{prop:2.5}, for these faithful integrals, we find that for any element $a\in A$ there is an element $b\in A$ so that $\varphi(ax)=\varphi(xb)$ for all $x$. The element $b$ is uniquely determined and we write it as $\sigma(a)$. We clearly have a linear map satisfying
\begin{equation*}
\varphi(x\sigma(aa'))=\varphi(aa'x)=\varphi(a'x\sigma(a))=\varphi(x\sigma(a)\sigma(a'))
\end{equation*}
so that $\sigma(aa')=\sigma(a)\sigma(a')$ by the faithfulness of $\varphi$. The faithfulness also implies that $\sigma$ is injective. And a similar argument to find $\sigma$ will also give that it is surjective. Hence, $\sigma$ is an automorphism of $A$.
\ssnl
Finally, if $a,b\in A$ we find
\begin{equation*}
\varphi(ab)=\varphi(b\sigma(a))=\varphi(\sigma(a)\sigma(b))=\varphi(\sigma(ab)).
\end{equation*}• 
Because $A$ is idempotent, it follows that $\varphi$ is invariant under $\sigma$.
\ebew

 We call it the {\it modular automorphism} for $\varphi$.
  \ssnl
 Observe that any faithful linear functional on a finite-dimensional algebra admits a modular automorphism.  In the infinite-dimensional case however, not all faithful linear functionals will admit a modular automorphism. So, in order to obtain the result above, we really need more than just a faithful linear functional.
 \ssnl
 In the literature, it is called the Nakayama automorphism (or rather its inverse), see e.g.\ \cite{Na}. The terms KMS property and modular automorphisms  find their origin in the theory of operator algebras and their usage in mathematical physics, see e.g.\ Chapter VIII in \cite{Ta}.
 \nl
A next result gives the relation between  left integrals.

\prop\label{prop:2.7}
Assume that $\varphi$ and $\varphi_1$ are left integrals and that $\varphi$ is  faithful. Then there is an element  $y$ in $A_s$ so that $\varphi_1(x)=\varphi(xy)$ for all $x$. 
\eprop

\bew
If we combine item i) of Proposition \ref{prop:2.5} with a single faithful right integral and these two left integrals, we find that for all $a$, there is an element $b$ so that $\varphi_1(xa)=\varphi(xb)$ for all $x$. We must have that $b=ay$ for a right multiplier $y$ of $A$. Then $\varphi_1(x)=\varphi(xy)$ for all $x$. In a similar way, we find a left multiplier $y'$ so that $\varphi_1(x)=\varphi(y'x)$ for all $x$. In particular $\varphi(xy)=\varphi(y'x)$ for all $x$.
\ssnl
If we replace $x$ by $ab$ and we get
\begin{equation*}
\varphi(ay\sigma(b))=\varphi(bay)=\varphi(y'ba)=\varphi(a\sigma(y'b)).
\end{equation*}
Replace $a$ by $a'a$ and use that this is true for all $a'$. It follows that $ay\sigma(b)=a\sigma(y'b)$ for all $a,b$. This is implies that $y\in M(A)$ and that $y=\sigma(y')$.
\ssnl
Now we apply Proposition \ref{prop:2.4} for $\varphi_1$ and for $\varphi$. We find 
\begin{equation*}
(\iota\ot\varphi)((1\ot a)\Delta(by))
=S((\iota\ot\varphi)(\Delta(a)(1\ot by)))=
(\iota\ot\varphi)((1\ot a)\Delta(b)(1\ot y)).
\end{equation*}
By the faithfulness of $\varphi$ we get $\Delta(by)=\Delta(b)(1\ot y)$ for all $b$. This implies that $y$ is in $A_s$.
\ebew

There is also a converse result. If $\varphi$ is any left integral and if $y\in A_s$, then $\varphi(\,\cdot\, y)$ as well as $\varphi(y\,\cdot\,)$, are again left invariant linear functionals.  
\snl
As a consequence, this will imply that $\sigma$ will leave $A_s$ globally invariant. 
\snl
Also, we find that there is a distinguished invertible element $y\in A_s$ satisfying $\varphi(S^2(x))=\varphi(xy)$ for all $x$.
\nl
Using completely similar methods, we can prove the following formula for a right integral in terms of a faithful left integral.

\prop\label{prop:2.8}
Let $\varphi$ be a faithful left integral and $\psi$ any right integral. Then there is an element $\delta\in M(A)$ such that $\psi(x)=\varphi(x\delta)$ for all $x$. If $\psi$ is also faithful, $\delta$ is invertible in $M(A)$.
\eprop

In many examples (as we will see e.g.\ in Section \ref{s:examples}) a single faithful integral exists. However, there are known examples where enough integrals exist (in the sense of our Definition \ref{defin:2.7} below), but not a single faithful one. See e.g.\ Proposition 2.5 in \cite{I-K}. 
\snl
Still, it is expected that the results above, proven in the case where a single faithful integral exists, will also have similar counterparts in the general situation.
\snl
We will discuss this further in Section \ref{s:conclusions}.

\nl
\bf The dual space $\widehat A$ \rm
\nl
We will now define the subspace $\widehat A$ of the dual space $A'$. Before we do that, we show in the following example that we can not expect that integrals on a weak multiplier Hopf algebra are automatically faithful (as in the case of multiplier Hopf algebras). This will be important for the upcoming definition of $\widehat A$ in this setting. 

\voorb
Take two Hopf algebras $(B,\Delta_B)$ and $(C,\Delta_C)$. Let $A$ be the direct sum of the algebras $B$ and $C$. So elements in $A$ are pairs $(b,c)$ with $b\in B$ and $c\in C$, with pointwise operations. We consider $B$ and $C$ as sitting in $A$ via the homomorphisms $b\mapsto (b,0)$ and $c\mapsto (0,c)$. The element $(1,0)$ is denoted by $e$ and then $(0,1)$ is $1-e$.
\ssnl
Define a coproduct $\Delta$ on $A$ by the requirement that it coincides with $\Delta_B$ on $B$ and with $\Delta_C$ on $C$.
The pair $(A,\Delta)$ is no longer a Hopf algebra, but a weak Hopf algebra. Indeed we get
\begin{equation*}
	\Delta(1)=\Delta(e)+\Delta(1-e)= e\ot e + (1-e)\ot (1-e).
\end{equation*}
We see that $\Delta(1)$ is strictly smaller than $1\ot 1$ in $A\ot A$.
\ssnl
If one of the components has integrals while the other has not, then we will have an example of a weak Hopf algebra with integrals, but obviously not enough integrals. More precisely, there will be no faithful integral.
\evoorb

So indeed, as an immediate consequence we see from this example that it will not be possible to show that an integral is automatically faithful as in the case of multiplier Hopf algebras. The example also makes clear that it will not be sufficient to assume the existence of integrals. We need {\it enough} integrals in the sense of the following definition.

\defin\label{defin:2.7}
We say that a weak multiplier Hopf algebra has a {\it faithful set of integrals} if the following two conditions are satisfied. Given an element $x\in A$, we must have $x=0$ if $\varphi(xa)=0$ for all left integrals $\varphi$ and elements $a\in A$. Similarly also if $\varphi(ax)=0$ for all left integrals $\varphi$ and elements $a\in A$, then $x=0$.
\edefin

If there is only one left integral (up to a scalar), as in the case of multiplier Hopf algebras, this condition is the same as faithfulness of this left integral. On the other hand, if there is a faithful left integral, then the condition above is fulfilled. So, what we essentially do with this definition is generalizing the requirement that there is a faithful left integral (or simply that there is an integral in the case of multiplier Hopf algebras).
\snl
Remark that we can formulate the condition also in terms of right integrals by applying the antipode.
\nl
Now, we are ready to give and discuss the definition of the {\it dual space} $\widehat A$ as follows.

\defin\label{defin:2.8}
Assume that $(A,\Delta)$ is a regular weak multiplier Hopf algebra with a faithful set of integrals. Then we define $\widehat A$ as the space of linear functionals on $A$ spanned by elements of the form $\varphi(\,\cdot\,a)$ where $\varphi$ is a left integral and $a\in A$.
\edefin

The choice of the representation of the elements in $\widehat A$ is not important as we will now prove, using the results of Proposition \ref{prop:2.5}.

\prop\label{prop:2.9}
For any left integral $\varphi$ and right integral $\psi$ and for all elements $a\in A$, we have that the linear functionals of the form
$\varphi(\,\cdot\,a)$, $\varphi(a\,\cdot\,)$, $\psi(\,\cdot\,a)$ en $\psi(a\,\cdot\,)$ are all in $\widehat A$. Moreover $\widehat A$ is the linear span of functionals of one of these four types, where $\varphi$, respectively $\psi$ runs over all left, respectively right integrals and $a$ over all elements in $A$.
\eprop

\bew 
By the assumption the set of left integrals is faithful as in Definition \ref{defin:2.7}. If we apply the antipode, we also find that the set of right integrals is faithful. 
\snl
We claim that all elements in $A$ can be obtained as a linear span of elements of one of the eight forms as we find them in the formulation of Proposition \ref{prop:2.5}. Let us show this for elements of the first form, namely 
\begin{equation*}
	(\iota\ot\varphi)((\iota\ot S)(\Delta(p))(1\ot q))
\end{equation*}•
where $p,q\in A$. If this is not be true, there exists a non-zero linear functional $\omega$ on $A$ that is $0$ on all such elements. For such a functional $\omega$, we have that $\varphi(S(r)q)=0$ with $r=(\omega\ot\iota)\Delta(p)$ for all $p,q\in A$ and for any left integral $\varphi$. Replace $q$ by a product $qq'$ and use the faithfulness of the set of left integrals to obtain that then $S(r)q=0$ for all $p,q$. Using the bijectivity of the antipode we obtain  that $(\omega\ot\iota)((1\ot q)\Delta(p))=0$, again for all $p,q$. From the fullness of the coproduct, it follows that $\omega=0$ and this is a contradiction. A similar argument is possible for all the other cases as we know that also the set of right integrals is faithful.
\snl
Then, using the results of Proposition \ref{prop:2.5}, it not only follows that elements of the four forms in the proposition all belong to $\widehat A$, but also that $\widehat A$ is the span of such elements for each of these forms.
\ebew

Observe that we used an argument that we already had in the beginning of the previous subsection, just before Proposition \ref{prop:2.7a}.
\snl
In what follows, we will use the following terminology.

\defin\label{defin:2.10}
Let $(A,\Delta)$ be a regular weak multiplier Hopf algebra. If there exists a faithful set of left integrals, we call it an {\it algebraic quantum groupoid} (or a {\it weak algebraic quantum group}).
\edefin

Recall that we also use the term {\it algebraic quantum group} for a {\it regular multiplier Hopf algebra with integrals}. The terminology is therefore consistent. 
\snl
Remark that the term {\it algebraic quantum groupoid} is also used in \cite{T2} within the context of multiplier Hopf algebroids. The notion is different from the one we introduce above in Definition \ref{defin:2.10}.
\snl
In Section 3, where we treat duality, we will use the above results to show that the adjoints of the product and coproduct of $A$ make $\widehat A$ again into a regular weak multiplier Hopf algebra with a faithful set of integrals. In other words, the dual of an algebraic quantum groupoid is again an algebraic quantum groupoid (see Theorem \ref{thm:3.18} in Section \ref{s:duality}). 

\nl
\bf Relation with the condition of Timmerman in \cite{T2}\rm
\nl
We finish this section with a look at the extra condition imposed in the work of Timmermann about the integrals on regular weak multiplier Hopf algebras when discussing the duality for these (see Definition 5.1.3 in \cite{T2}).
\snl
As before, we assume that $(A,\Delta)$ is a regular weak multiplier Hopf algebra. We have the canonical idempotent $E$ and the source and target algebras $\varepsilon_s(A)$ and $\varepsilon_t(A)$.
\snl
First we prove two lemmas. The first one is about a faithful linear functional on $A$. There is no need for it to be an integral.

\lem
Let $f$ be a faithful linear functional on $A$.  Then any element in $\varepsilon_t(A)$ is of any of the forms
\begin{equation*}
(f\ot\iota)((a\ot 1)E), \tussen (f\ot\iota)(E(a\ot 1)),
\end{equation*}
for some $a\in A$. On the other hand, any element in $\varepsilon_s(A)$ is of any of the forms 
\begin{equation*}
(\iota\ot f)((1\ot a)E), \tussen (\iota\ot f)(E(1\ot a)),
\end{equation*}
for some $a\in A$. 
\elem

\bew
We will prove it for elements of the form $(f\ot\iota)((a\ot 1)E)$. The proof in the three other cases is completely similar.
\ssnl
First observe that $(a\ot 1)E$ is an element of $A\ot \varepsilon_t(A)$. Indeed, we have
\begin{equation*}
(pq\ot 1)E=\sum_{(q)}pq_{(1)}\ot S^{-1}(q_{(3)})q_{(2)}=\sum_{(q)}pq_{(1)}\ot S^{-1}(\varepsilon_s(q_{(2)}))
\end{equation*}
and this belongs to $A\ot \varepsilon_t(A)$ for all $p,q$. As $A$ is idempotent, we get $(a\ot 1)E\in A\ot \varepsilon_t(A)$ for all $a$ in $A$. See e.g.\ Formula (1.6) in \cite{VD-W5}. Similar results like this are also found already in the appendix of \cite{VD-W4}. 
\ssnl
Let $\omega$ be a linear functional on $\varepsilon_t(A)$ and assume that it is $0$ on $(f\ot\iota)((a\ot 1)E)$ for all $a\in A$. Then $f(ax)=0$ for all $a$ where $x=(\iota\ot\omega)E$. This element $x$ belongs to $M(A)$. Then also $f(aa'x)=0$ for all $a,a'\in A$ and because $f$ is assumed to be faithful, we must have $a'x=0$ for all $a'\in A$. In other words $(\iota\ot\omega)((a'\ot 1)E)=0$ for all $a'\in A$. By the fullness of $E$, we have that $\omega=0$ on all of $\varepsilon_t(A)$. This proves that any element in $\varepsilon_t(A)$ is of the form $(f\ot\iota)((a\ot 1)E)$ for some $a\in A$.
\ebew

In the next lemma, we have some kind of a converse. Now we need to work with integrals.
\lem
Let $\varphi$ be a left integral. It will be faithful if
any element in $\varepsilon_t(A)$ is of the form $(\varphi\ot\iota)((a\ot 1)E)$ for some $a$ in $A$ and also of the form 
$(\varphi\ot\iota)(E(a\ot 1))$ for some $a$.
Similarly, let $\psi$ be a right integral. It will be faithful if any element in $\varepsilon_s(A)$ is of the form $(\iota\ot\psi)((1\ot a)E)$ for some $a$ and of the form $(\iota\ot\psi)(E(1\ot a))$ for some $a$.
\elem

\bew
i) Let $\varphi$ be a left integral. Let $b\in A$ and assume that $\varphi(ab)=0$ for all $a\in A$. Then we have
\begin{equation*}
(\iota\ot\varphi)(\Delta(a)(1\ot b))=0
\end{equation*}
for all $a\in A$. By the result in Proposition \ref{prop:2.4}, we also find that 
\begin{equation*}
(\iota\ot\varphi)((1\ot a)\Delta(b))=0
\end{equation*}
for all $a\in A$. Now we apply $\Delta$, use the Sweedler notation for the coproduct on $b$ and replace $a$ by $aS(b_{(2)})$. This give that 
\begin{equation*}
\sum_{(b)} b_{(1)}\varphi(aS( b_{(2)}) b_{(3)})=0
\end{equation*}
for all $a$. Now, as in the the previous proof, we have 
\begin{equation*}
\sum_{(b)} b_{(1)}\ot S( b_{(2)}) b_{(3)}=(\iota\ot S)((b\ot 1)E).
\end{equation*}
It follows that $by=0$ where $y=(\iota\ot\varphi)((1\ot a)(\iota\ot S)E)$. Now we use that $(S\ot S)E=\zeta E$ where $\zeta$ is the flip. We get $y=S^{-1}(x)$ where $x=(\varphi\ot\iota)((a\ot 1)E)$. By assumption any element in $\varepsilon_t(A)$ is of this form and because $S$ is a bijection from $\varepsilon_t(A)$ to $\varepsilon_s(A)$ we finally get $by=0$ for all $y\in \varepsilon_s(A)$. We know that this implies that $b=0$. See e.g.\ Proposition 2.10 and Lemma 2.12 in \cite{VD-W5}. 
\ssnl
ii) Next assume that $\varphi(ba)=0$ for all $a\in A$. In a completely similar way, it will follow that $b=0$ if we know that any element in $\varepsilon_s(A)$ is of the form $(\varphi\ot\iota)(E(a\ot 1))$ for some $a$ in $A$.
\ssnl
iii) This proves the first statement of the lemma. The result for a right integral can be proven in a similar way from the two other conditions. But in fact, it follows from the result for a left integral by using the fact that the antipode converts a right integral into a left integral.
\ebew

Before we continue, remark that this argument is precisely the one that is used to show that a non-zero left integral on regular multiplier Hopf algebras is automatically faithful. See e.g.\ Proposition 3.4 in \cite{VD2}. Indeed, in that case $E=1$ and $\varepsilon_s(A)$ are just the scalar multiples of the identity and the condition of the lemma is automatically fulfilled as soon as $\varphi$ is non-zero.
\snl
Now we combine the two results and  we obtain the following criteria for faithfulness of integrals, when they exist.

\prop
Let $(A,\Delta)$ be a regular weak multiplier Hopf algebra. Assume that $\varphi$ is a left integral on $A$. Then it is faithful if and only if 
\begin{align*}
\{(\varphi\ot\iota)((a\ot 1)E)\mid a\in A\}&=\varepsilon_t(A) \\
\{(\varphi\ot\iota)(E(a\ot 1))\mid a\in A\}&=\varepsilon_t(A). 
\end{align*}
On the other hand, let $\psi$ be a right integral on $A$. Then it is faithful if and only if
\begin{align*}
\{(\iota\ot\psi)((1\ot a)E)\mid a\in A\}&=\varepsilon_s(A) \\
\{(\iota\ot\psi)(E(1\ot a))\mid a\in A\}&=\varepsilon_s(A) \\
\end{align*}
\eprop
\vskip -10pt
So we see that the conditions imposed in Definition 5.1.3 in \cite{T2} are equivalent with requiring the existence of a single faithful left integral. Remark again that the second statement about right integrals is equivalent with the first statement on left integrals because the antipode converts one to the other, flips $E$ and converts the set $\varepsilon_s(A)$ to $\varepsilon_t(A)$.
\snl

It is also not hard to verify the following, more general result.

\prop
Assume that we have a set $\Cal S$ of left integrals. This will be a faithful set if and only if $\varepsilon_t(A)$ is spanned by elements of the form
\begin{equation*}
(\varphi\ot\iota)((a\ot 1)E)
\end{equation*}
where $\varphi\in\Cal S$ and $a\in A$, as well as by elements of the form
\begin{equation*}
(\varphi\ot\iota)(E(a\ot 1))
\end{equation*}
where again $\varphi\in\Cal S$ and $a\in A$.
\eprop

This property is indeed more general as it can happen that a weak multiplier Hopf algebras has a faithful set of integrals, but not a single faithful one.
\nl\nl

\section{\hspace{-17pt}. Duality} \label{s:duality} 

Let $(A,\Delta)$ be a regular weak multiplier Hopf algebra. Assume that there is a faithful set of integrals (as defined in Definition \ref{defin:2.7}). So, using the terminology as introduced in Section \ref{s:integrals}, we assume that $(A,\Delta)$ is an {\it algebraic quantum groupoid} (cf.\ Definition \ref{defin:2.10}). 
\snl
Recall the definition of the dual $\widehat A$. It is the space of linear functionals spanned by elements of the form $\varphi(\,\cdot\,a)$ where $\varphi$ is a left integral and $a$ an element of $A$ (see Definition \ref{defin:2.8}). Remember that we also can take elements of the form $\varphi(a\,\cdot\,)$ and that left integrals can be replaced by right integrals (see Proposition \ref{prop:2.9}). The dual space $\widehat A$ separates points of $A$ by definition (because the set of left integrals is assumed to be faithful). 
\snl
In this section, we will show that the space $\widehat A$ carries a product (dual to the coproduct on $A$) and a coproduct $\widehat \Delta$ (dual to the product on $A$) so that the pair $(\widehat A,\widehat \Delta)$ is again a regular weak multiplier Hopf algebra. As expected, the antipode $\widehat S$ on $\widehat A$ will simply be the adjoint of $S$.  We will give a formula for the dual canonical idempotent $\widehat E$. We will also show that this dual weak multiplier Hopf algebra again has a faithful set of (left) integrals and we will provide formulas for those. So, the pair $(\widehat A,\widehat \Delta)$ is again an algebraic quantum groupoid, called the dual of $(A,\Delta)$. The dual of $(\widehat A,\widehat \Delta)$ is canonically isomorphic with the original weak multiplier Hopf algebra $(A,\Delta)$.
\snl
In the case where there is a single faithful (left) integral,  $\widehat A$ has a single faithful integral as well. In that case it should be possible to give explicit formulas for the dual objects $\widehat \delta$, $\widehat \sigma$, ... in terms of the objects of the original weak multiplier Hopf algebras. These formulas are expected to be very similar (if not completely the same) as in the  case of algebraic quantum groups (and algebraic quantum hypergroups). See e.g.\ Proposition 4.1 in \cite{D-VD}. We refer to Section \ref{s:conclusions} for more comments on this case.
\snl
The case where there is a single faithful integral is also studied in Section 5 of \cite{T2} as we have mentioned before. It is shown that also the dual has a single faithful integral in that case. We include an argument in this section as well. See Theorem \ref{stel:2.20} below. 
\nl
\bf The dual algebra $\widehat A$ \rm
\nl
We begin with the definition of the product on $\widehat A$.

\prop\label{prop:3.1}
For any pair $\omega,\omega'$ of elements in $\widehat A$, we can define the product $\omega\omega'\in\widehat A$ by
\begin{equation*}
	(\omega\omega')(x)=(\omega\ot\omega')(\Delta(x))
\end{equation*}
where $x\in A$. If $\omega'=\varphi(a\,\cdot\,)$, where $\varphi$ is a left integral and $a\in A$, then $\omega\omega'=\varphi(b\,\cdot\,)$ where now 
\begin{equation*}
b=((\omega\circ S)\ot\iota)\Delta(a).
\end{equation*}
\eprop

\bew
The proof is standard. First we must argue that this product is well-defined as a linear functional on $A$. Take $\omega'=\sum_i\varphi_i(\,\cdot\,a_i)$ where $\varphi_i$ are left integrals and $a_i\in A$.  Then we see that
\begin{equation*}
	(\omega\ot\omega')\Delta(x)=\sum_i(\omega\ot\varphi_i)(\Delta(x)(1\ot a_i))
\end{equation*}
and this is defined for all $x$ because for all $a_i$ we have that $\Delta(x)(1\ot a_i)$ belongs to $A\ot A$. One must also argue that $\omega\omega'$ is well-defined in the sense that it does not depend on the choice of the representation of $\omega'$. For doing this, one assumes that $\omega$ is chosen as $\varphi(\,\cdot\,a)$ for a left integral $\varphi$ and an element $a\in A$.
\ssnl
Next one has to show that the resulting element $\omega\omega'$ again belongs to $\widehat A$. To see this, take $\omega'=\varphi(a\,\cdot\,)$ where $\varphi$ is a left integral and $a\in A$. Using the formula in Proposition \ref{prop:2.4}, we find
\begin{equation*}
	(\omega\ot\omega')\Delta(x)=(\omega\ot\varphi)((1\ot a)\Delta(x))
=((\omega\circ S)\ot\varphi)(\Delta(a)(1\ot x))
\end{equation*}
and we see that $\omega\omega'=\varphi(b\,\cdot\,)$ with $b=((\omega\circ S)\ot\iota)\Delta(a)$. This element $b$ belongs to $A$ as also $\omega\in\widehat A$. 
\ssnl
This completes the proof.
\ebew

Also the proof of the following result is not so difficult.

\prop\label{prop:3.2}
The product defined on $\widehat A$ in the previous proposition makes $\widehat A$ into an associative algebra with a non-degenerate product. The algebra is also idempotent.
\eprop

\bew
The {\it associativity} of the product on $\widehat A$ follows of course from the coassociativity of the coproduct $\Delta$ on $A$. In order to do this in a careful way, using the notion of coassociativity as defined in Definition 1.1 of \cite{VD-W4}, one can use the various representations of elements in $\widehat A$ as obtained in Proposition \ref{prop:2.9} of Section \ref{s:integrals}. 
\ssnl
Next we show that the product is {\it non-degenerate}. Assume first that $\omega\in \widehat A$ and that $\omega\omega'=0$ for all $\omega'\in\widehat A$. This implies that 
$\omega'((\omega\ot\iota) \Delta(x))=0$
for all $x\in A$ and for all  $\omega'\in \widehat A$. 
Recall that $(\omega\ot\iota) \Delta(x)$ is in $A$ for all $x\in A$ because $\omega\in \widehat A$. 
As $\widehat A$ is separating points of $A$, it follows that $(\omega\ot\iota) \Delta(x)=0$ for all $x\in A$. If we apply the counit, we find that $\omega(x)=0$ for all $x$ and this means that $\omega=0$.
Similarly, we will get $\omega=0$ if $\omega\omega'=0$ for all $\omega'\in\widehat A$.
\ssnl
Finally we prove that the algebra is {\it idempotent}. To show this, let $\omega\in \widehat A$, $a\in A$ and let $\varphi$ be any left integral. Denote $\omega'=\varphi(a\,\cdot\,)$. As in the proof of Proposition \ref{prop:3.1}, we find $\omega\omega'=\varphi(b\,\cdot\,)$ with $b=((\omega\circ S)\ot\iota)\Delta(a)$. Now we claim that any element in $A$ is a linear combination of elements $b$ of this form. This will imply that any element in $\widehat A$ can be written as a sum of products of elements in $\widehat A$. 
\ssnl
To prove the claim, assume that this is not the case. Then there exists a linear functional $\rho$ on $A$ so that 
$((\omega\circ S)\ot \rho)\Delta(a)=0$ for all $\omega \in \widehat A$ and all $a\in A$. We now proceed as above (when proving that the product is non-degenerate), but we need to be a little more careful. Replace $\omega$ by $\omega(\,\cdot\,b)$ with $b\in A$ and use that $\widehat A$ separates points. This will give that $(S\ot \rho)(\Delta(a))(b\ot 1)=0$ for all $a,b$. This means that $(\iota\ot \rho)\Delta(a)=0$ in $M(A)$ for all $a$. Using that the coproduct is full will give that $\rho=0$. This completes the proof.
\ebew

Remark that in the second part of the proof above, in stead of using the counit, we can also proceed by using that the coproduct is full (as we had to do in the last part of the proof where we could not apply the counit). 
\snl
Also observe that for proving that $\widehat A$ is idempotent, we do not need that the original algebra $A$ is idempotent.
\snl

The proofs of the previous results on the dual algebra are completely the same as in the case of multiplier Hopf algebras with integrals (see \cite{VD2}) and algebraic quantum hypergroups (see \cite{D-VD}). The construction of the coproduct is however, just as in these two cases, a little more subtle. We treat it carefully and give more details.
\nl
\bf The pairing $\langle A,\widehat A\rangle$ \rm
\nl
Before we investigate the coproduct on $\widehat A$, let us first have a closer look at the non-degenerate pairing between the algebras $A$ and $\widehat A$. As an application, we will get a characterization of $M(\widehat A)$ that is similar as for algebraic quantum (hyper)groups, (see e.g.\ Proposition 3.4 in \cite{D-VD}). 
\snl
Further in this section, we will also use the {\it pairing notation}. So we will write $\langle a,\omega \rangle$ for $\omega(a)$ when $a\in A$ and $\omega\in\widehat A$. In this case, we will  use $B$ for $\widehat A$ and $a,a',\dots$ for elements in $A$ and $b,b',\dots$ for elements in $B$. 
\snl
Left and right multiplication in the algebras on either side, induce right and left actions as we see in the following proposition.

\prop\label{prop:3.3}
There exist left and right actions of $A$ on $B$, as well as right and left actions of $B$ on $A$ given by the following formulas:
\begin{align*}
	\langle a'a,b\rangle&=\langle a',a\tr b\rangle \qquad\qquad \langle a,bb'\rangle=\langle a\tl b,b'\rangle\\
	\langle aa',b\rangle&=\langle a', b\tl a\rangle \qquad\qquad \langle a,b'b\rangle=\langle b\tr a,b'\rangle
\end{align*}
where $a,a'\in A$ and $b,b'\in B$. In the two cases, the left and right actions commute. Moreover, all these four actions are unital and non-degenerate.
\eprop

\bew
First consider the {\it left action of $A$ on $B$}. Assume that $\omega=\varphi(\,\cdot\,a')$. Then we have $a\tr \omega=\omega(\,\cdot\,a)=\varphi(\,\cdot\,aa')$ and this will again belong to $\widehat A$. We also see that $A\tr B=A$ because $A^2=A$. So the action is unital. Finally, to show that the action is non-degenerate, assume that $\omega\in \widehat A$ and that $\omega(\,\cdot\,a)=0$ for all $a$. Again because $A^2=A$, it follows that $\omega=0$. The arguments for the right action of $A$ on $B$ are completely similar. It is clear that the left action commutes with the right action.
\ssnl
Next consider the {\it right action of $B$ on $A$}. We have that 
\begin{equation*} 
a\tl \omega=\textstyle\sum_{(a)}\omega(a_{(1)}) a_{(2)}
\end{equation*}
where $\omega\in \widehat A$. We know that this belongs to $A$. We get indeed a right action of $\widehat A$ on $A$. Similarly, we get the left action. They clearly commute.
\ssnl
To show that also this right action is unital, we have to prove that any element in $A$ is the span of such elements. Suppose that this is not the case. Then, there exists a non-zero linear functional $\rho$ on $A$ such that $(\omega\ot \rho)\Delta(a)=0$ for all $a\in A$ and all $\omega\in \widehat A$. Then, as in the proof of Proposition \ref{prop:3.2}, it follows that $\rho=0$ and this gives a contradiction. Finally, also this action is non-degenerate because the pairing is non-degenerate and $B$ is also idempotent (see Proposition \ref{prop:3.2}). Similarly for the left action.
\ebew 

Observe that we know the existence of local units for $A$ (see Proposition 4.8 in \cite{VD-W4}). Then any unital action of $A$ is non-degenerate. However, at this moment, we do not yet know that also $\widehat A$ has local units. This will only follow when we have shown that it is also a regular weak multiplier Hopf algebra. This means that, at this level, we have to show that these actions are unital {\it and} non-degenerate.
\snl
Next, we will use the results of Proposition \ref{prop:3.3} to extend the pairing from $A \times B$ to $A\times M(B)$ and to $M(A)\times B$. The two cases are treated in a different way, just as above. Remark e.g.\ that extending the pairing to $M(A)\times B$ is a special case of extending a {\it reduced} linear functional on $A$ to $M(A)$. It is known that this can be done by using that $A$ has local units. The other case is different.
\snl
First we have the easier case.

\prop\label{prop:3.4}
There is a unique extension of the pairing of $A$ with $B$ to a pairing of $M(A)$ with $B$, satisfying 
\begin{equation*}
	\langle m, b\tl a\rangle = \langle am, b \rangle
	\tussenen
	\langle m, a\tr b\rangle = \langle ma,b\rangle
\end{equation*}
whenever $a\in A$, $b\in B$  and $m\in M(A)$.
\eprop

\bew
The proof is standard. Start with an element $m\in M(A)$ and $b\in B$. We use e.g.\ the fact that the right action $\tl$ of $A$ on $B$ is unital to write $b=\sum_i b_i\tl a_i$. We try to define 
\begin{equation*}
	\langle m, b \rangle=\sum_i\langle a_im,b_i\rangle,
\end{equation*}
but we first have to show that this is well-defined. For this assume that $\sum_i b_i\tl a_i=0$. Take $e$ in $A$ so that $a_ie=a_i$ for all $i$. Then
\begin{equation*}
	\sum_i\langle a_im,b_i\rangle=\sum_i\langle a_iem,b_i\rangle=\sum_i\langle em,b_i\tl a_i\rangle=0.
\end{equation*}
This proves that we can define the pairing, using the first formula.
 To prove that the second formula also holds, one replaces $b$ by $\sum b_i\tl a_i$, uses the first formula and again that the right action is unital. Uniqueness is obvious.
\ebew

We now look at the other extension. Remark again that we used local units above and so, we have to find another way to prove this case.

\prop\label{prop:3.5}
The pairing on $A\times B$ can be extended to $A\times M(B)$ in such a way that 
\begin{equation*}
	\langle a,bm\rangle=\langle a\tl b,m\rangle 
	\tussenen
	\langle a,mb\rangle=\langle b\tr a,m\rangle
\end{equation*}
whenever $a\in A$, $b\in B$ and $m\in M(B)$.
\eprop

\bew
We start as in the previous proof. Now take $m\in M(B)$, $a\in A$, use the fact that the right action $\tl$ of $B$ on $A$ is unital and write $a=\sum_i a_i\tl b_i$. We again want to define 
\begin{equation*}
	\langle a,m\rangle=\sum_i \langle a_i,b_im\rangle.
\end{equation*}
So, we assume that $\sum_i a_i\tl b_i=0$. Then we have, for all $b\in B$,
\begin{align*}
 	\sum_i\langle a_i\tl (b_im),b\rangle 
	&=\sum_i\langle a_i,b_imb\rangle \\
	&=\sum_i\langle a_i\tl b_i,mb\rangle=0.
\end{align*}
This implies that $\sum_i a_i\tl (b_im)=0$. If we apply the counit $\varepsilon$ of $A$ we find precisely that 
\begin{equation*}
	\sum_i \langle a_i,b_im\rangle=0.
\end{equation*}
This shows that we can define the extension and that it will satisfy the first formula. An argument, similar as in the previous proof will give also the second formula. And again, uniqueness is trivial.
\ebew 

So, instead of using the property of having local units, we use the counit. Of course, the counit is the identity in $M(B)$ and so, in a sense, we use that we already have the extension of the pairing to elements of the form $(a,1)$ in $A\times M(B)$. Indeed, we are using that $\varepsilon(a\tl b)=\langle a, b\rangle$. But this can be justified by replacing $b$ by an element of the form $a'\tr b$. 
\snl
We will use this result later, in Proposition \ref{prop:3.9} below, to define the counit $\widehat\varepsilon$ on $\widehat A$ by $\widehat\varepsilon(\omega)=\omega(1)$ for $\omega\in\widehat A$.
\snl
Recall that in the general case (just as for multiplier Hopf algebras), the pairing cannot be extended further to $M(A)\times M(\widehat A)$.
\nl
Before we continue with the study of the coproduct on $\widehat A$, we apply the above results to get the following characterization of elements in $M(\widehat A)$. We have the same characterization as for algebraic quantum (hyper-)groups, see e.g.\ Proposition 3.4 in \cite{D-VD}.

\prop
Let $\omega$ be a linear functional in $A'$. Then there is a multiplier $m$ in $M(\widehat A)$ so that $\omega(a)=\langle a,m\rangle$ if and only if
\begin{equation*}
	(\omega\ot\iota)\Delta(a)
	\tussenen
	(\iota\ot\omega)\Delta(a)
\end{equation*}
are in $A$ for all $a\in A$.
\eprop

\bew
i) First let $m\in M(\widehat A)$ and define $\omega$ in $A'$ by $\omega(a)=\langle a,m\rangle$. We use the extended pairing as obtained in the previous proposition. Now we use e.g.\ that the left action of $\widehat A$ on $A$ is unital. Then we can write any element $a\in A$ as $\sum_ib_i\tr a_i$. This will give 
\begin{equation*}
	(\iota\ot\omega)\Delta(a)=\sum_i (mb_i)\tr a_i
\end{equation*}
and this belongs to $A$. Similarly, if $a=\sum_i b_i\tl a_i$, we find
\begin{equation*}
(\omega\ot\iota)\Delta(a)=\sum_i a_i \tl (b_im)
\end{equation*}• and we have that also $(\omega\ot\iota)\Delta(a)$ belongs to $A$.
\snl
ii) Conversely, assume that we have an element $\omega$ in $A'$ satisfying the conditions in the formulation of the proposition. Then we can define a multiplier $m$ of $\widehat A$ by
\begin{equation*}
	(m\omega')(a)=(\omega\ot\omega')\Delta(a)
	\tussenen
	(\omega'm)(a)=(\omega'\ot\omega)\Delta(a)
\end{equation*}
where $\omega'\in \widehat A$. Arguments as used earlier (see the proof of Proposition \ref{prop:3.1}) will be needed in order to have that $m\omega'$ and $\omega'm$ again belong to $\widehat A$. The rest of the argument is more or less trivial.
\ebew

The proof is essentially the same as for an algebraic quantum group and even an algebraic quantum hypergroup. See e.g.\ Proposition 3.4 in \cite{D-VD}. 
\nl
\bf The coproduct $\widehat\Delta$ on $\widehat A$\rm
\nl
Now, we construct the coproduct. We do it by first defining the canonical maps $\widehat{T_1}$ and $\widehat {T_2}$ in the following proposition. As before, we use the pairing and we write $\langle a,\omega \rangle$ for $\omega(a)$ when $a\in A$ and $\omega\in\widehat A$ (or even in $A'$). Similarly, we use this notation when $a\in A\ot A$ and $\omega$ in $\widehat A\ot \widehat A$ (or even $\omega\in (A\ot A)')$.
\prop\label{prop:3.7}
There is a regular coproduct $\widehat \Delta$ on $\widehat A$ so that the associated canonical maps $\widehat{T_1}$ and $\widehat {T_2}$, satisfy 
\begin{align}
	\langle x\ot y,\widehat{T_1}(\omega\ot\omega') \rangle 
							&=\langle T_2(x\ot y),\omega\ot\omega' \rangle \label{eqn:3.1}\\
       \langle x\ot y,\widehat{T_2}(\omega\ot\omega') \rangle 
							&=\langle T_1(x\ot y),\omega\ot\omega' \rangle\label{eqn:3.2}
\end{align}
for all $x,y\in A$ and $\omega,\omega'\in \widehat A$.
\eprop

\bew
i) The formulas (\ref{eqn:3.1}) and (\ref{eqn:3.2}) in the formulation of the proposition define $\widehat{T_1}$ and $\widehat {T_2}$ as linear maps from $\widehat A\ot \widehat A$ to the linear dual space $(A\ot A)'$ of $A\ot A$. We claim that these maps have range in $\widehat A\ot \widehat A$. 
\snl
To show this for $\widehat{T_1}$, we take 
$\omega'=\varphi(a\,\cdot\,)$ where $\varphi$ is a left integral on $A$ and where $a\in A$. Then we have for all $x\in A$ and $y\in A$ that
\begin{align*}
	\langle T_2(x\ot y),\omega\ot\omega' \rangle 
		&= (\omega\ot\varphi)((x\ot a)\Delta(y)) \\
		&= (\omega\ot\varphi)((x\ot 1)((S\ot\iota)(\Delta(a)(1\ot y)))). 
\end{align*}
We see that 
\begin{equation}
\widehat{T_1}(\omega\ot\omega')=\sum_{(a)}\omega(\,\cdot\,S(a_{(1)}))\ot \varphi(a_{(2)}\,\cdot\,)\label{eqn:3.4b}
\end{equation}
and so $\widehat{T_1}(\omega\ot\omega') \in \widehat A\ot \widehat A$
because $((S\ot\iota)\Delta(a))(a'\ot 1) \in A\ot A$ when $a,a'\in A$. To prove that $\widehat {T_2}$ maps $\widehat A\ot \widehat A$, we use appropriate representations of the elements $\omega,\omega'$ of $\widehat A$, using right integrals.
\snl
ii) Next, we show that there is a coproduct $\widehat\Delta$ on $\widehat A$ so that $\widehat{T_1}$ and $\widehat{T_2}$ are the associated canonical maps. For this we need that 
\begin{equation}
	(\omega\ot 1)(\widehat{T_1}(\omega'\ot\omega''))=(\widehat{T_2}(\omega\ot\omega'))(1\ot\omega'') \label{eqn:3.3}
\end{equation}
for all $\omega$, $\omega'$ and $\omega''$ in $\widehat A$. Remark that this equation makes sense in $\widehat A\ot \widehat A$. And indeed, it will imply that there is a map 
\begin{equation*}
	\widehat\Delta:\widehat A\to M(\widehat A\ot\widehat A) 
\end{equation*}
so that $\widehat{T_1}(\omega'\ot\omega'')=\widehat\Delta(\omega')(1\ot \omega'')$ and $\widehat{T_2}(\omega\ot\omega')=(\omega\ot 1)\widehat\Delta(\omega')$ for all $\omega$, $\omega'$ and $\omega''$ in $\widehat A$. 
\snl
To prove that (\ref{eqn:3.3}) holds, one can  verify that
\begin{align*}
	\langle x\ot y, (\omega\ot 1)(\widehat{T_1}(\omega'\ot \omega''))\rangle 
	&= \langle \Delta(x)\ot y, (\iota\ot \widehat{T_1})(\omega\ot\omega'\ot \omega'')\rangle \\
	&= \langle (\iota\ot T_2)(\Delta(x)\ot y), \omega\ot\omega'\ot \omega'' \rangle \\
	&= \langle (\Delta(x)\ot 1)(1\ot \Delta(y)), \omega\ot\omega'\ot \omega'' \rangle
\end{align*}
for all $x,y\in A$ and $\omega$, $\omega'$ and $\omega''$ in $\widehat A$. Similarly also
\begin{equation*}
	\langle x\ot y , (\widehat{T_2}(\omega\ot\omega'))(1\ot \omega''\rangle
 = \langle (\Delta(x)\ot 1)(1\ot \Delta(y)), \omega\ot\omega'\ot \omega'' \rangle
\end{equation*}
for all such elements. Some care is needed, but remark that in all these equations, the legs of $\Delta(x)$ and the legs of $\Delta(y)$ are well covered by using the appropriate forms of the elements in $\widehat A$.
\snl
iii) That this coproduct is regular, is proven as in i), using again well-chosen representations of the elements in $\widehat A$.
\snl
iv) Next we argue that $\widehat\Delta$ is a homomorphism. We will do this with the pairing notation. For notational convenience, we will again use $B$ for $\widehat A$ and denote the coproduct $\widehat\Delta$ on $B$ simply by $\Delta$. We will also use the Sweedler notation for the two coproducts. The reader can verify that the necessary coverings exist (possibly also induced by the pairing and {\it repeated} coverings).
\snl
Take $a,a'$ in $A$ and $b,b',c$ in $B$. Then we have
\begin{align*}
	\langle a\ot a',\Delta(bb')(1\ot c)\rangle
	&=\langle (a\ot 1)\Delta(a'),bb'\ot c\rangle \\
	&=\sum_{(a')}\langle aa'_{(1)},bb'\rangle\langle a'_{(2)},c\rangle \\
	&=\sum_{(a')}\langle \Delta(aa'_{(1)}),b\ot b'\rangle\langle a'_{(2)},c\rangle \\
	&=\sum_{(a')}\langle \Delta(aa'_{(1)})\ot  a'_{(2)},b\ot b'\ot c\rangle.
\end{align*}
We have used that the map $\widehat T_1$ is the adjoint of $T_2$. Also remark that $c$ will provide the necessary coverings.
\snl
On the other hand we find
\begin{align*}
	\langle a\ot a',\Delta(b)\Delta(b')(1\ot c)\rangle
	&=\sum_{(b),(b')}\langle a\ot a',b_{(1)}b'_{(1)}\ot b_{(2)}b'_{(2)}c)\rangle\\
	&=\sum_{(b),(b')}\langle \Delta(a)\ot a',b_{(1)}\ot b'_{(1)}\ot b_{(2)}b'_{(2)}c)\rangle\\
	&=\sum_{(b')}\langle (\Delta(a)\ot 1)\Delta_{13}(a'),b\ot b'_{(1)}\ot b'_{(2)}c)\rangle\\
	&=\langle (\Delta(a)\ot 1)\Delta^{(2)}(a'),b\ot b'\ot c)\rangle.
\end{align*}
Again we have used that $\widehat T_1$ is the adjoint of $T_2$, now two times. And again, all the time, $c$ will take care of the necessary coverings.
\ssnl
We see that both expressions are the same and this proves that $\Delta$ is a homomorphism on $B$.
\snl
v)) It remains to be shown that $\widehat\Delta$ is coassociative. For this we need to show that the maps $\iota\ot\widehat{T_1}$ and $\widehat{T_2}\ot \iota$ commute on $\widehat A\ot\widehat A$. This will follow by duality if we can show that the maps $\iota\ot T_2$ and $T_1\ot \iota$ commute on $A\ot A$. And this in turn should follow from the {\it associativity} of the product in $A$. Indeed, we get for all $x,y,z$ in $A$ that
\begin{align*}
	 (T_1\ot\iota)(\iota\ot T_2)(x\ot y\ot z)
	&=(T_1\ot\iota)(x\ot(y\ot 1)\Delta(z)))\\
	&=(\Delta(x)\ot 1)(1\ot y\ot 1)(1\ot\Delta(z))
\end{align*}•
and similarly for $(\iota\ot T_2)(T_1\ot\iota)$.
\ebew

Also this proof is essentially the same as for regular multiplier Hopf algebras with integrals. Another similar treatment is found in Section 3 of \cite{D-VD} in the case of algebraic quantum hypergroups. We essentially use only the formulas from Proposition  \ref{prop:2.4}. They are also valid for algebraic quantum hypergroups. That is why this works in all these cases in the same way.
\snl
Observe in the last part of the proof, that the dual maps $\iota\ot\widehat{T_1}$ and $\widehat{T_2}\ot \iota$ commute, is equivalent with the fact that the maps $\iota\ot T_2$ and $T_1\ot \iota$ commute on $A\ot A$. This last fact is different from the property that the maps $\iota\ot T_1$ and $T_2\ot \iota$ commute on $A\ot A$. The first commutation rule is equivalent with the {\it associativity of the product} in $A$, whereas the latter with the {\it coassociativity of the coproduct} on $A$.
\snl
We finish this item on the coproduct $\widehat\Delta$ by proving the following expected  formula. Observe that we need the extension of the pairing from $(A\ot A)\times (\widehat A\ot\widehat A)$ to $(A\ot A)\times M(\widehat A\ot\widehat A)$. This extension can be obtained as in Proposition \ref{prop:3.5}. We again use the pairing conventions.

\prop\label{prop:3.10}
For all $a,a'$ in $A$ and $b\in B$, we have $\langle aa',b\rangle=\langle a\ot a',\Delta(b)\rangle$.
\snl\bf Proof\rm:
Take $a,a' \in A$ and $b,b'\in B$. Then we have
\begin{align*}
	\langle a\ot (b'\tr a'),\Delta(b)\rangle
	&= \langle a\ot a',\Delta(b)(1\ot b'\rangle \\
	&= \langle (a\ot 1)\Delta(a'),b\ot b'\rangle \\
	&= \langle a(b'\tr a'),b\rangle.
\end{align*}
Now the result follows because the left action of $B$ on $A$ is unital.
\eprop

\nl
\bf Existence of a counit and fullness of the coproduct \rm
\nl
Next, we show that there exists a counit $\widehat\varepsilon$ and that the coproduct $\widehat\Delta$ is full (so that hence the counit is unique). This is what we show first.

\prop\label{prop:3.8}
The coproduct $\widehat\Delta$ on $\widehat A$, as obtained in Proposition \ref{prop:3.7}, is full.
\eprop

\bew
We have seen in the proof of the Proposition \ref{prop:3.7} that 
\begin{equation*}
	\widehat\Delta(\omega)(1\ot \omega')=\sum_{(a)}\omega(\,\cdot\,S(a_{(1)}))\ot\varphi(a_{(2)}\,\cdot\,)
\end{equation*}
if $\omega,\omega'\in \widehat A$ and if $\omega'$ has the form $\varphi(a\,\cdot\,)$ for $a\in A$ and a left integral $\varphi$ on $A$. See Equation (\ref{eqn:3.4b}). If $\gamma$ is the linear functional on $\widehat A$, given by the evaluation in a point $c$ of $A$, we find that
\begin{equation*}
	(\iota\ot\gamma(\,\cdot\,\omega'))\widehat\Delta(\omega)
	=\sum_{(a)}\omega(\,\cdot\,S(a_{(1)})))\varphi(a_{(2)}c).
\end{equation*}
So, in order to show that the 'left leg' of $\widehat\Delta$ is all of $\widehat A$, we must have that any element in $A$ is a linear combination of elements of the form
\begin{equation*}
	\sum_{(a)}S(a_{(1)})\varphi(a_{(2)}c)
\end{equation*}
where $a,c\in A$ and where $\varphi$ is a left integral on $A$. Results of this type been obtained before, see e.g.\ the proof of Proposition \ref{prop:2.9}.
\ssnl
A similar argument works  for the right leg.
\ebew

\prop\label{prop:3.9}
There is a unique counit $\widehat\varepsilon$ on $\widehat A$ given by $\widehat\varepsilon(\omega)=\omega(1)$.
\eprop

\bew 
First remark that we can use Proposition \ref{prop:3.4} (with $m=1$ in $M(A)$) in order to define $\omega(1)$ for $\omega\in \widehat A$.  It is also possible to consider $\varphi(a)$ when $\omega=\varphi(\,\cdot\, a)$ but then one would need to show that this is a well-defined map from $\widehat A$ to $\mathbb C$. The proof of this can be given by using local units for $A$, just as we needed for the proof of Proposition \ref{prop:3.4}.
\ssnl
By definition, we will have
$\widehat\varepsilon(\omega(\,\cdot\,a))=\omega(a)$
as well as
$\widehat\varepsilon(\omega(a\,\cdot\,))=\omega(a)$
for all $a\in A$ and $\omega\in\widehat A$.
\ssnl
Consider $\omega$ and $\omega'$ in $\widehat A$. We know by definition that 
\begin{equation*}
	(\widehat\Delta(\omega)(1\ot \omega'))(a\ot a')=(\omega\ot\omega')((a\ot 1)\Delta(a')).
\end{equation*}
Therefore we find that 
\begin{equation*}
	(\widehat\varepsilon\ot\iota)(\widehat\Delta(\omega)(1\ot \omega'))=(\omega\ot\omega')\circ\Delta=\omega\omega'.
\end{equation*}
We need Proposition \ref{prop:3.10} here. 
Similarly we will get 
\begin{equation*}
	(\iota\ot\widehat\varepsilon)((\omega\ot 1)\widehat\Delta(\omega'))=\omega\omega'.
\end{equation*}
This proves that $\widehat\varepsilon$ is a counit. And because we know already that the coproduct is full, it is the unique coproduct on the pair $(\widehat A,\widehat\Delta)$.
\ebew

\nl
\bf The antipode $\widehat S$ on $\widehat A$ \rm
\nl
We now look for the appropriate generalized inverses $\widehat R_1$ and $\widehat R_2$ of $\widehat T_1$ and $\widehat T_2$ respectively. Because we have that $\widehat T_1$ and $\widehat T_2$ are adjoint to $T_2$ and $T_1$, we expect that $\widehat R_1$ and $\widehat R_2$ are adjoint to $R_2$ and $R_1$ respectively. Of course, we also expect that the antipode $\widehat S$ of the pair $(\widehat A, \widehat\Delta)$ will be the adjoint of $S$.
\nl
First we consider the adjoint of $S$.

\prop\label{prop:3.11}
There is a map $\widehat S:\widehat A \to \widehat A$, satisfying
$$\langle a, \widehat S(b)\rangle = \langle S(a),b \rangle$$
for all $a\in A$ and $b\in \widehat A$. This map is also bijective.
\eprop

\bew
The above formula can be used to define $\widehat S(b)\in A'$ for all $b\in A'$. We just have to argue that $\widehat S(b)\in \widehat A$ for all $b\in \widehat A$.
\snl
So take e.g.\ $b=\varphi(\,\cdot\,c)$ where $c\in A$ and $\varphi$ is a left integral on $A$. Then
\begin{equation*}
	\langle a, \widehat S(b) \rangle
	=\varphi(S(a)c)=(\varphi\circ S)(S^{-1}(c)a)
\end{equation*}
for all $a$ so that $\widehat S(b)=(\varphi\circ S)(S^{-1}(c)\,\cdot\,)$. We know that $\varphi\circ S$ is a right integral. From Proposition \ref{prop:2.9} we find that $\widehat S(b)$ again belongs to $\widehat A$.
\snl
We also see from the above argument that $\widehat S$ is a bijective map from $\widehat A$ to itself.
\ebew

As the map $S$ is both an anti-algebra and an anti-coalgebra map, the same will be true for its adjoint (because the product and the coproduct are dual to each other). This is standard. 
\snl
Remark that at this moment, we still have not shown that the pair $(\widehat A, \widehat\Delta)$ is again a weak multiplier Hopf algebra and so we cannot yet argue that the map $\widehat S$, defined above, is actually the antipode of $(\widehat A, \widehat\Delta)$. This will become obvious later.
\snl
A similar remark applies for the maps $\widehat R_1$ and $\widehat R_2$ we consider next.

\prop\label{prop:3.12}
There exist maps $\widehat R_1$ and $\widehat R_2$ from $\widehat A\ot\widehat A$ to itself, satisfying
\begin{align}
	\langle a\ot a',\widehat R_1(b\ot b')\rangle &= \langle R_2(a\ot a'),b\ot b' \rangle  \label{eqn:3.4}\\
	\langle a\ot a',\widehat R_2(b\ot b')\rangle  &= \langle R_1(a\ot a'),b\ot b' \rangle  \label{eqn:3.5}
\end{align}
for all $a,a'\in A$ and $b,b'\in \widehat A$. Moreover we have
\begin{align}
	\widehat R_1(b\ot b')&=\sum_{(b)}b_{(1)} \ot \widehat S(b_{(2)}) b' \label{eqn:3.6} \\
	\widehat R_2(b\ot b')&=\sum_{(b')}b\widehat S(b'_{(1)}) \ot b'_{(2)} \label{eqn:3.7} 
\end{align}
for all $b,b'$. The maps $\widehat R_1$ and $\widehat R_2$ are generalized inverses of $\widehat T_1$ and $\widehat T_2$ respectively.
\eprop

\bew
Again we can use the formulas (\ref{eqn:3.4}) and (\ref{eqn:3.5}) above to define $\widehat R_1(b\ot b')$ and $\widehat R_2(b\ot b')$ in $(A\ot A)'$ for $b,b'\in A'$. It can be shown as in the proof of Proposition \ref{prop:3.7} that actually these elements belong to $\widehat A\ot \widehat A$ if $b,b'\in \widehat A$. In fact, it is also possible to use that $\widehat S$ is an anti-automorphism of the algebra $\widehat A$, combined with the fact that the canonical maps on $\widehat A\ot \widehat A$ have range in $\widehat A\ot \widehat A$.
This takes care of the first statement of the proposition.
\ssnl
Let us now prove Equation (\ref{eqn:3.6}). The proof of (\ref{eqn:3.7}) is completely similar. 
\ssnl
Take $a,a'\in A$ and $b,b'\in \widehat A$.
If we consider the map $\widehat R_1$ as in the formula (\ref{eqn:3.6}), we get
\begin{align*}
	\langle a\ot a',\widehat R_1(b\ot b')\rangle
	&= \langle a\ot \Delta(a'), (\iota\ot \widehat S\ot\iota)(\widehat \Delta(b)\ot b')\rangle\\
	&= \langle a\ot (S\ot \iota)\Delta(a'), \widehat \Delta(b)\ot b'\rangle \\
	&= \langle (a\ot 1)(S\ot \iota)\Delta(a'), b\ot b'\rangle \\
	&= \langle R_2(a\ot a'), b\ot b'\rangle.
\end{align*}
This proves the equality (\ref{eqn:3.6}).
\ssnl
Finally, it easily follows that $\widehat R_1$ is a generalized inverse of $\widehat T_1$ because they are adjoints of $R_2$ and $T_2$ that themselves are generalized inverses of each other. Similarly  $\widehat R_2$ is a generalized inverse of $\widehat T_2$. This completes the proof.
\ebew

One can check that the necessary coverings can be found for the formulas above.
\snl
Remark that in the case of a multiplier Hopf algebra, this is a lot simpler. Then we just have that these maps are the inverses of the canonical maps. We know that they are bijective in that case.
\snl
Also observe, as we already have the counit $\widehat\varepsilon$ available, that it follows that the map $\widehat S$ satisfies the expected formulas
\begin{equation*}
	\sum_{(b)}b_{(1)}\widehat S(b_{(2)})b_{(3)}=b
	\tussenen
	\sum_{(b)}\widehat S(b_{(1)})b_{(2)}\widehat S(b_{(3)})=\widehat S(b)
\end{equation*}
for all $b$ in $\widehat A$.

\nl
\bf The canonical idempotent $\widehat E$ \rm
\nl
For convenience in what follows, we will not always write  $\widehat \Delta$, but often use also $\Delta$ for the dual coproduct. Similarly for the dual canonical maps and there generalized inverses. At the same time, we also keep using $B$ for $\widehat A$ and $b,b',\dots$ for elements in the dual. 
\snl
It is not hard to guess what the canonical idempotent $\widehat E$ should be. Indeed, eventually it has to be $\Delta(1)$ and so we expect that $\widehat E$, as an element in the dual space $(A\ot A)'$, should satisfy
\begin{equation*}
	\langle a\ot a',\widehat E\rangle = \langle aa', 1 \rangle =\varepsilon(aa')
\end{equation*}
\snl
We will now use this formula to {\it define} $\widehat E$ and prove its properties. In particular, we show that this is indeed the canonical idempotent for the dual $(\widehat A,\widehat \Delta)$.

\prop\label{prop:3.13}
There is a multiplier $\widehat E$ in $M(B\ot B)$ defined by
\begin{equation*}
	\langle a\ot a',\widehat E\rangle=\varepsilon(aa')
\end{equation*}
for all $a,a'\in A$. It satisfies
\begin{equation*}
	\widehat E(b\ot b')=T_1R_1(b\ot b')
	\tussenen
	(b\ot b')\widehat E=T_2R_2(b\ot b')
\end{equation*}
for all $b,b'\in B$. In particular, it is an idempotent multiplier determining the ranges of the  canonical maps as in the definition of a weak multiplier Hopf algebra (Definition 1.14 in \cite{VD-W4}).
\eprop

\bew
We first define $\widehat E$ in $(A\ot A)'$ using the above formula. Next take two elements $b,b'$ in $B$. Then, using the product in $(A\ot A)'$, we get
\begin{equation}
	\langle a\ot a', \widehat E(b\ot b')\rangle
	=\sum_{(a)(a')}\varepsilon(a_{(1)}a'_{(1)}) \langle a_{(2)},b\rangle \langle a'_{(2)},b' \rangle.\label{eqn:3.8}
\end{equation}
Now we have
\begin{equation*}
\sum_{(a')} \Delta(a)(a'_{(1)}\ot 1)\ot a'_{(2)}=
\sum_{(a')} \Delta(aa'_{(1)})(1\ot S(a'_{(2)}))\ot a'_{(3)}
\end{equation*}
and if we apply $\varepsilon$ on the first factor, we find
\begin{equation*}
\sum_{(a),(a')} \varepsilon(a_{(1)}a'_{(1)})a_{(2)}\ot a'_{(2)}=
\sum_{(a')} aa'_{(1)}S(a'_{(2)})\ot a'_{(3)}=R_2T_2(a\ot a').
\end{equation*}
If we use this in the formula (\ref{eqn:3.8}), we find
\begin{equation*}
	\langle a\ot a', \widehat E(b\ot b')\rangle =\langle R_2T_2(a\ot a'),b\ot b'\rangle
\end{equation*}
for all $a,a'$ and $b,b'$. This shows that $\widehat E(b\ot b')=T_1R_1(b\ot b')$. In a similar way, we find that also $(b\ot b')\widehat E=T_2R_2(b\ot b')$. 
\ebew

Next we show that the coproduct has the right behavior on the legs of the canonical multiplier $\widehat E$.

\prop\label{prop:3.14}
We have
\begin{align*}
	(\Delta\ot\iota)\widehat E &= (\widehat E\ot 1)(1\ot \widehat E)\\
	(\Delta\ot\iota)\widehat E &= (1\ot \widehat E)(\widehat E \ot 1).
\end{align*}
\eprop

\bew
To prove the first equality, multiply from the left with $\Delta(b)\ot b'$ where $b,b'\in B$. Then we see that it is equivalent with showing that
\begin{equation*}
	(\Delta\ot\iota)T_2R_2=(\iota\ot T_2R_2)(\Delta\ot\iota).
\end{equation*}
This equation is adjoint to the equation
\begin{equation*}
	R_1T_1(aa'\ot a'')=(a\ot 1)R_1T_1(a'\ot a'')
\end{equation*}
for $a,a',a''\in A$. And because $R_1T_1(a\ot a')=(a\ot 1)F_1(1\ot a')$ for all $a,a'\in A$ where $F_1=(\iota\ot S)E$,  (see Propositions 4.5 and 4.7 in \cite{VD-W4}), this equation is true.
\snl
To prove the other equality $(\Delta\ot\iota)\widehat E=(1\ot \widehat E)(\widehat E\ot 1)$, we multiply from the right with $\Delta(b)\ot b'$. Then the formula is equivalent with
\begin{equation*}
	(\Delta\ot\iota)T_1R_1=(\iota\ot T_1R_1)(\Delta\ot\iota).
\end{equation*}•
This in turn is adjoint to the equation
\begin{equation*}
 R_2T_2(aa'\ot a'')=(a\ot 1)R_2T_2(a'\ot a'')
\end{equation*}
for all $a,a',a''$. And this is again true because $R_2T_2(a\ot a')=(a\ot 1)F_2(1\ot a')$ with $F_2=(S\ot\iota)E$, see again Propositions 4.5 and 4.7 in \cite{VD-W4}.
\ebew

Now we are ready for the main result.

\stel\label{thm:3.15}
The dual pair $(\widehat A,\widehat \Delta)$ is again a regular weak multiplier Hopf algebra.
\estel

\bew
We will use Theorem  2.9 of \cite{VD-W4}.
\ssnl
i) We have shown in Proposition \ref{prop:3.2} that $\widehat A$ is a non-degenerate, idempotent algebra. We also have proven in Proposition \ref{prop:3.7} and \ref{prop:3.8} that there is a regular and full  coproduct $\widehat \Delta$ on $\widehat A$. And finally, in Proposition \ref{prop:3.9},  that this coproduct admits a counit.
\ssnl
ii) We have defined the linear map $\widehat S$ from $\widehat A$ to itself (in Proposition \ref{prop:3.11}) and the associated maps $\widehat R_1$ and $\widehat R_2$ with the formulas (\ref{eqn:3.6}) and (\ref{eqn:3.7}) in Proposition \ref{prop:3.12}. They are canonical inverses of $\widehat T_1$ and $\widehat T_2$ respectively. It follows that the map $\widehat S$ satisfies the requirements in item i) and the equalities (2.5) in Theorem 2.9 of \cite{VD-W4}. 
\ssnl
iii) Finally, we have found the dual canonical multiplier $\widehat E$ satisfying the necessary formulas as proven in Proposition \ref{prop:3.13} and \ref{prop:3.14}.
\vskip 5pt
Hence all conditions of Theorem 2.9 in \cite{VD-W4} are fulfilled and we have that the pair $(\widehat A,\widehat \Delta)$ is a weak multiplier Hopf algebra. We know that $\widehat S$ is its antipode and because this is a bijective map from $\widehat A$ to itself (see \ref{prop:3.11}), the weak multiplier Hopf algebra $(\widehat A,\widehat \Delta)$ is regular.
\ssnl
This completes the proof.
\ebew
As we mentioned already, in the case of a single faithful integral, this result is obtained in \cite{T2} as a consequence of the duality for regular multiplier Hopf algebroids with integrals and the relation of multiplier Hopf algebroids with weak multiplier Hopf algebras as obtained in \cite{T-VD1}. We refer to Section \ref{s:conclusions} for a discussion on this approach.
\nl
\bf Existence of integrals on the dual \rm
\nl
Next we will show that the dual weak multiplier Hopf algebra also has a faithful set of integrals.

\prop\label{prop:3.16}
For every element $a\in A$, there exists a right invariant functional $\psi_a$ on $B$ such that
\begin{equation}
	\psi_a(\omega)=\varphi(a\varepsilon_s(c))\qquad\text{ when } \qquad \omega=\varphi(\,\cdot\, c)\label{eqn:3.9} 
\end{equation}
for any left integral $\varphi$ on $A$ and any element $c \in A$. 
\eprop

\bew
Fix an element $a\in A$.
\ssnl
i) We first show that we can define $\psi_a$ such that (\ref{eqn:3.9}) holds. To do this, assume that we have a finite number of left integrals $\varphi_i$ and elements $c_i\in A$ so that
\begin{equation*}
	\sum_i \varphi_i(\,\cdot\,c_i)=0.
\end{equation*}
This implies that
\begin{equation*}
	\sum_i(\iota\ot\varphi_i)(\Delta(x)(1\ot c_i))=0
\end{equation*}
for all $x\in A$. Apply $S$ and use one the formulas in Proposition \ref{prop:2.4}, to obtain that also
\begin{equation*}
	\sum_i((\iota\ot\varphi_i)((1\ot x)\Delta(c_i))=0
\end{equation*}
for all $x$. Again apply $S$ and multiply with the element $a$. Using the Sweedler notation we find
\begin{equation*}
	\sum_{i,(c_i)}  aS(c_{i(1)}) \ot \varphi_i(\,\cdot\, c_{i(2)})=0.
\end{equation*}
If we apply the evaluation map $p\ot\omega\mapsto \omega(p)$, to this equation we get
\begin{equation*}
	\sum_{i,(c_i)}\varphi_i(aS(c_{i(1)}) c_{i(2)})
	=\sum_i\varphi_i(a\varepsilon_s(c_i))=0.
\end{equation*}
It follows that we can define a linear map $\psi_a$ on $B$ as in (\ref{eqn:3.9}) above.
\ssnl
ii) Next we show that the functional $\psi_a$ is right invariant on $B$. To do this, take a left integral $\varphi$ on $A$ and two elements $c,c'\in A$. Let $\omega=\varphi(\,\cdot\,c)$ as before. Then
\begin{equation*}
	\psi_a(\omega(\,\cdot\,c'))
	=\varphi(a\varepsilon_s(c'c))
	=\varphi(a\varepsilon_s(\varepsilon_s(c')c))
\end{equation*}
where we have used that $\varepsilon_s(c'c))=\varepsilon_s(\varepsilon_s(c')c)$. This equality follows e.g.\ from 
\begin{equation*}
\Delta(\varepsilon_s(c')c)=(1\ot \varepsilon_s(c'))\Delta(c).
\end{equation*}

We see that
\begin{equation}
	\psi_a(\omega(\,\cdot\,c'))=\psi_a(\omega(\,\cdot\,\varepsilon_s(c')).\label{eqn:2.11}
\end{equation}
We will now show that this implies right invariance of $\psi_a$. 
\ssnl
We claim that 
\begin{equation*}
\omega(\,\cdot\,\varepsilon_s(c'))=\langle (\omega\ot\iota)\widehat F_1,\,\cdot\,\ot c'\rangle.
\end{equation*}
Then equation (\ref{eqn:2.11}) above means that
\begin{equation*}
(\psi_a\ot\iota)\widehat\Delta(\omega)=(\psi_a\ot\iota)((\omega\ot\iota)\widehat F_1)
\end{equation*}
and this will show that $\psi_a$ is right invariant.
\ssnl
To prove the claim, take $p\in A$ and $c'$ of the form $\omega'\tr q$ for some $\omega'\in\widehat A$ and $q\in A$. Then
\begin{align*}
\langle (\omega\ot\iota)\widehat F_1,p\ot c' \rangle
&=\langle (\omega\ot\iota)\widehat F_1,p\ot \omega'\tr q \rangle\\
&=\langle (\omega\ot\iota)\widehat F_1(1\ot \omega',p\ot q \rangle\\
&=\langle R_1T_1(\omega\ot\omega',p\ot q \rangle\\
&=\langle \omega\ot\omega',T_2R_2 (p\ot q) \rangle\\
&=\textstyle\sum_{(q)} \omega(p\varepsilon_s(q_{(1)})\omega'(q_{(2)})
=\omega(p\varepsilon_s(c')).
\end{align*}
This proves the claim and completes the proof of right invariance of $\psi_a$
\ebew

In the next section we have an example where the coproduct on the dual actually maps into the algebraic tensor product and then, this property is more easy to obtain. We refer to \cite{VD8} where this is explained in greater detail.
\snl
A more formal argument would be as follows.
\snl
Using the Sweedler notation, the pairing notation and writing $b$ for $\omega$ we can rewrite Equation (\ref{eqn:2.11}) and obtain
\begin{align*}
	\sum_{(b)}\psi_a(b_{(1})\langle c',b_{(2)}\rangle
	&=\sum_{(b)}\psi_a(b_{(1})\langle \varepsilon_s(c'),b_{(2)}\rangle \\
	&=\sum_{(b)}\psi_a(b_{(1})\langle c',\varepsilon_s(b_{(2)})\rangle.
\end{align*}
We need here that the source map is  {\it self-dual}. Then we obtain
\begin{equation*}
	\sum_{(b)}\psi_a(b_{(1})b_{(2)}
	=\sum_{(b)}\psi_a(b_{(1)})\varepsilon_s(b_{(2)}),
\end{equation*}
showing also the right invariance of $\psi_a$. This argument makes sense e.g.\ when the source maps have range in the algebra itself, and not just in the multiplier algebra. But that would more or less imply that the algebras are finite-dimensional. To make this argument precise seems not so easy.
\snl
For completeness, we show that the source and target maps are self-dual maps below.

\prop\label{prop:2.17}
For all $a\in A$ and $b\in B$ we have
\begin{align*}
\langle a,\varepsilon_s(b) \rangle&=\langle\varepsilon_s(a),b\rangle\\
\langle a,\varepsilon_t(b) \rangle&=\langle\varepsilon_t(a),b\rangle
\end{align*}
\eprop

The formulas make sense using the extended pairing, discussed earlier in this section. 
\snl
It is also easy to give a formal proof. Indeed, given $a\in A$ and $b\in B$, we have
\begin{align*}
\langle a,\varepsilon_s(b) \rangle
&=\sum_{(b)} \langle a, S(b_{(1)}) b_{(2)}\rangle\\
&=\sum_{(a),(b)}\langle a_{(1)}\ot a_{(2)}, S(b_{(1)})\ot b_{(2)}\rangle\\
&=\sum_{(a),(b)}\langle S(a_{(1)})\ot a_{(2)}, b_{(1)}\ot b_{(2)}\rangle\\
&=\sum_{(a)}\langle S(a_{(1)})a_{(2)}, b\rangle
=\langle\varepsilon_s(a),b\rangle.
\end{align*}
We see that $\varepsilon_s$ is a self-dual map. Similarly we have
\begin{equation*}
\langle a,\varepsilon_t(b) \rangle=\langle\varepsilon_t(a),b\rangle
\end{equation*}
for all $a\in A$ and $b\in B$.
\snl
However, it is not trivial to see that all this makes sense. The problem is to get coverings for the intermediate steps. We give an argument in the proof below.
\snl

\vspace{-0.3cm}\begin{itemize}\item[ ] \bf Proof \rm (of Proposition \ref{prop:2.17}): 
Start with elements $p,a\in A$ and $q,b\in B$. We know that
\begin{equation}
\langle T_2R_2(p\ot a),q\ot b \rangle=\langle p\ot a,R_1T_1(q\ot b)\rangle.\label{eqn:2.12}
\end{equation}
Using the Sweedler notation, we find for the left hand side
\begin{align*}
\langle T_2R_2(p\ot a),q\ot b \rangle
&=\sum_{(a)} \langle pS(a_{(1)})a_{(2)} \ot a_{(3)},q\ot b \rangle \\
&=\langle p\varepsilon_s(b\tr a),q \rangle.
\end{align*}
For the right hand side we get
\begin{align*}
\langle p\ot a,R_1T_1(q\ot b)\rangle
&=\sum_{(q)}\langle p\ot a, q_{(1)}\ot S(q_{(2)})q_{(3)}b\rangle\\
&=\langle a,  \varepsilon_s(q\tl p)b \rangle.
\end{align*}
So we can rewrite Equation (\ref{eqn:2.12}) as
\begin{equation*}
\langle p\varepsilon_s(b\tr a),q \rangle=\langle a,  \varepsilon_s(q\tl p)b \rangle
\end{equation*}
and this means 
\begin{equation*}
\langle \varepsilon_s(b\tr a), q\tl p \rangle=\langle b\tr a,\varepsilon(q\tl p)\rangle.
\end{equation*}
The result now follows from the fact that these actions are unital.
\ebew

\snl
Next, we show that the set of right integrals is faithful.

\prop\label{prop:3.17}
The set of right integrals $\psi_a$ with $a\in A$, as defined in the previous proposition, is faithful (in the sense of Definition \ref{defin:2.7}).
\eprop

\bew
i) Take $\omega=\varphi(\,\cdot\,c)$  where $\varphi$ is a left integral on $A$ and $c\in A$. Let $\omega'\in \widehat A$. Then for all $p\in A$ we have
\begin{align*}
(\omega'\omega)(p)
&=\sum_{(p)}\omega'(p_{(1)})\varphi(p_{(2)}c)\\
&=\sum_{(c)}\omega'(S^{-1}(c_{(1)}))\varphi(pc_{(2)}).
\end{align*}
By the definition of $\psi_a$ we find $\psi_a(\omega'\omega)=\sum_{(c)}\omega'(S^{-1}(c_{(1)})\varphi(a\varepsilon_s(c_{(2)}))$.
\snl
ii) Now assume that $\psi_a(w'w)=0$ for all $a$ and all $\omega$. Because we assume a faithful set of left integrals on $A$, it follows that $\sum_{(c)}\omega'(S^{-1}(c_{(1)}))\varepsilon_s(c_{(2)})=0$ for all $c$. We now apply the distinguished linear functional defined on $\varepsilon_s(A)$ by   $\varepsilon_s(p))\mapsto\varepsilon(p)$ for all $p$. This will give $\omega'(S^{(-1)}(c))=0$ for all $c$. Hence $\omega'=0$.
\snl
iii) On the other hand assume $\omega=\sum_i \varphi_i(\,\cdot\,c_i)$ and $\psi_a(\omega'\omega)=0$ for all $\omega'$. Then
\begin{equation*}
\sum_{i,(c_i)} \varphi_i(a\varepsilon_s(c_{i(2)}))S^{-1}(c_{i(1)})=0
\end{equation*}
for all $a$. We find
\begin{equation*}
x\mapsto \sum_{i,(c_i)} \varphi_i(ax\varepsilon_s(c_{i(2)})) \ot c_{i(1)}
\end{equation*}
for all $a$. Finally apply the evaluation map and use that $c_i=\sum_{(c_i)}c_{i(1)}\varepsilon_s(c_{i(2)})$ to arrive at
\begin{equation*}
\sum_{i,(c_i)} \varphi_i(ac_i) 
\end{equation*}
for all $a$. This proves that $\omega=0$.

\ebew

Now, we can prove the second main result.

\stel\label{thm:3.18}
The dual $(\widehat A,\widehat \Delta)$ of the weak multiplier Hopf algebra $(A,\Delta)$ again has a faithful set of integrals. Then we can consider the dual of $(\widehat A,\widehat \Delta)$ and it is canonically isomorphic with the original pair $(A,\Delta)$.
\estel

\bew
Most of this has been shown. Indeed, in Theorem \ref{thm:3.15}, we have shown that $(\widehat A,\widehat \Delta)$ is again a regular weak multiplier Hopf algebra. On the other hand, in Proposition \ref{prop:3.17} we found that it has a faithful set of integrals. So, the only thing that is left is the existence of a bijection of the dual of $\widehat A$ to $A$.
\ssnl
If we look again at the proof of Proposition \ref{prop:3.17} above, we see e.g.\ that 
$\psi_a(\,\cdot\,\omega')$ is given by evaluation in the point $\sum_{(c')}S^{-1}(c'_{(1)})\varphi(a\varepsilon_s(c'_{(2)}))$ 
where $\omega'=\varphi(\,\cdot\,c')$. We claim that these elements belong to $A$ and that all elements of $A$ are of this form.
\ssnl
First, recall that 
\begin{equation}
\sum_{(c')}c'_{(1)}\ot \varepsilon_s(c'_{(2)})b=R_1T_1(c'\ot b)=(c'\ot 1)F_1(1\ot b)
\end{equation}
for any $b\in A$. If we replace $c'$ by $cc'$ we see that
\begin{equation}
(cc'\ot 1)F_1=\sum_{(c')}cc'_{(1)}\ot \varepsilon_s(c'_{(2)})
\end{equation}
and it follows that $\sum_{(c')}c'_{(1)}\ot \varepsilon_s(c'_{(2)})$ is an element in $A\ot M(A)$ and so 
\begin{equation*}
\sum_{(c')}c'_{(1)}\ot a\varepsilon_s(c'_{(2)})
\end{equation*}
 is in $A\ot A$. This proves the first claim.
\ssnl
To prove the second part of the claim, assume that $\rho$ is a linear functional that vanishes on all such elements. By the fact that we have a faithful set of integrals on $A$, this will give that
$(\rho (S^{-1}(\,\cdot\,))\ot\iota )((c'\ot 1)F_1)=0$ for all $c'$. This implies that $\rho=0$ because $F_1=(\iota\ot S)E$ and $E$ is full.
\ebew

\nl
\bf The case of a single faithful integral - Relation with \cite{T2} \rm
\nl
We now assume that there is a single faithful left integral $\varphi$ on $A$. In this case, we can prove the following result about integrals on the dual.

\prop\label{prop:2.19}
Let $p$ be any linear functional on the source algebra $\varepsilon_s(A)$. Define $\widehat\psi$ on $\widehat A$ by
\begin{equation*}
\widehat\psi(\omega)=p(\varepsilon_s(c)) 
\qquad\text{ when } \qquad 
 \omega=\varphi(\,\cdot\,c)
\end{equation*}
where $c\in A$. Then $\widehat\psi$ is a right integral on the dual $(\widehat A,\widehat\Delta)$.
\eprop

\bew
i) In this case, there is no need to argue that $\widehat\psi$ is well defined. By the faithfulness of $\varphi$ it follows that $c$ is determined by $\omega$. Moreover, any element in $\widehat A$ is of the form $\varphi(\,\cdot\,c)$ for some $c\in A$.
\ssnl
ii) To prove right invariance, take $c\in A$ and $\omega=\varphi(\,\cdot\,c)$. As in the prove of Proposition \ref{prop:3.16} we have, for all $a\in A$, 
\begin{equation*}
\sum_{(\omega)} \omega_{(1)} \omega_{(2)}(a)=\varphi(\,\cdot\,ac).
\end{equation*}
Hence 
\begin{equation*}
\sum_{(\omega)}\widehat\psi(\omega_{(1)}) \omega_{(2)}(a)=p(\varepsilon_s(ac)).
\end{equation*}
Invariance of $\widehat\psi$ again follows from $\varepsilon_s(ac)=\varepsilon_s(\varepsilon_s(a)c)$.
\ebew

Observe the difference with the formula in Proposition \ref{prop:3.16}. There we also use a function on the source algebra, but since we are working with several left integrals, we need a certain relation of the functionals for each left integral. Otherwise, it is not possible to show that it is well defined.
\snl
On the other hand, we now can take the distinguished linear functional on $\varepsilon_s(A)$ for $p$. Then we get $\widehat\psi(\omega)=\varepsilon(c)$ for $\omega=\varphi(\,\cdot\,c)$ because now $p(\varepsilon_s(a)=\varepsilon(a)$. This give the right integral from Corollary 5.17 in \cite{T2}. For this choice, one has that $\widehat\psi$ is again faithful on $\widehat A$.  Indeed, by a simple calculation, we find that, if $\omega=\varphi(\,\cdot\,c)$, then $\widehat\psi(\omega'\omega)=\omega'(S^{-1}(c))$ for all $\omega'$. If this is equal to $0$ for all $\omega$, then $\omega'(S^{-1}(c))=0$ for all $c$ so that $\omega'=0$. Similarly, if this is $0$ for all $\omega'$, then $c$ and hence $\omega=0$.
\snl
So we get the following result (Corollary 5.17 in \cite{T2}).

\stel\label{stel:2.20}
If $(A,\Delta)$ is a regular weak multiplier Hopf algebra with a single faithful integral, then the dual $(\widehat A,\widehat\Delta)$ is again a regular weak multiplier Hopf algebra with a single faithful integral.
\estel
\nl
\bf The source and target of the dual \rm
\nl
We have seen that the source and target maps are self-dual maps. See Proposition \ref{prop:2.17}. 
There is however another important relation between the source and target on the dual  and the source and target on the original weak multiplier Hopf algebra. 
\prop\label{prop:3.19}
There is an isomorphism $\gamma_s$ from the source algebra $\varepsilon_s(A)$ of $A$ to the target algebra $\varepsilon_t(\widehat A)$ of the dual, given by the formula
\begin{equation}
\langle ya,b\rangle=\langle a,b\,\gamma_s(y) \rangle\label{eqn:3.10}. 
\end{equation}
for all $a\in A$ and $b\in B$. Here $y\in \varepsilon_s(A)$. Similarly, there is an isomorphism $\gamma_t$ from the target algebra $\varepsilon_t(A)$ of $A$ to the source algebra $\varepsilon_s(\widehat A)$ of the dual, given by the formula
\begin{equation}
\langle ax,b \rangle = \langle a ,\gamma_t(x)b \rangle
\end{equation}
for all $a\in A$ and $b\in B$. Here $x\in \varepsilon_t(A)$.
\eprop
\bew
We use that $T_1R_1(a\ot a')=E(a\ot a')$ for $a,a'\in A$ and that $R_2T_2(b\ot b')=(b\ot 1)F_2(1\ot b')$ when $b,b'\in \widehat A$ where $F_2=(S\ot\iota)\widehat E$. We also know that the canonical map $T_2$ on $\widehat A\ot \widehat A$ is adjoint to the canonical map $T_1$ on $A\ot A$ and that the generalized inverse $R_2$ of $T_2$ for the dual is in turn the adjoint of the generalized inverse $R_1$ of the original canonical map $T_1$. 
\snl
Now take $a,a'\in A$ and $b,b'\in \widehat A$. Then
\begin{align*}
\langle E(a\ot a'),b\ot b' \rangle
&=\langle T_1R_1(a\ot a'),b\ot b'  \rangle \\
&=\langle a\ot a', R_2T_2(b\ot b')  \rangle \\
&=\langle a\ot a', (b\ot 1)F_2(1\ot b')  \rangle. 
\end{align*}
We know that $E(1\ot a')$ belongs to $\varepsilon_s(A)\ot A$ and consequently, if we pair with an element $b'$ in the second factor, we get an element of $\varepsilon_s(A)$. From the fact that the pairing is non-degenerate, we obtain that all elements in $\varepsilon_s(A)$ can be obtained as a linear span of such elements. A similar dual result is true for the right hand side. 
\snl
Then, for all $y\in \varepsilon_s(A)$ there is a $y_1$ in $\varepsilon_s(\widehat A)$ so that 
\begin{equation*}
\langle ya,b \rangle = \langle a,bS(y_1) \rangle
\end{equation*}
for all $a\in A$ and $b\in B$. It is clear that $y_1$ is uniquely determined by $y$. Hence we get a linear map $\gamma_s:\varepsilon_s(A)\to \varepsilon_t(\widehat A)$ by putting $\gamma_s(y)=S(y_1)$. It is also clear that $\gamma_s$ is a homomorphism and that it is injective. A similar argument as before will give that it is also surjective. Hence we get an isomorphism $\gamma_s:\varepsilon_s(A)\to \varepsilon_t(\widehat A)$ determined by the formula (\ref{eqn:3.10}).
\snl
The other result is proven in a completely similar way. We now use that $$T_2R_2(a\ot a')=(a\ot a')E$$ and that $R_1T_1(b\ot b')=(b\ot 1)F_1(1\ot b')$ where $F_1=(\iota\ot S)E$.
\ebew

We can obtain partial results like above by using the following argument.
\snl
Take $y\in \varepsilon_s(A)$.
For $a\in A$  and $b,b'\in \widehat A$ we have on the one hand
\begin{align*}
\langle ya,bb' \rangle 
&= \langle \Delta(ya),b\ot b' \rangle \\
&= \langle (1\ot y)\Delta(a),b\ot b' \rangle \\
&=\langle \Delta(a),b\ot y\tr b' \rangle \\
&=\langle a, b(y\tr b' \rangle.
\end{align*}

On the other hand, this is equal to $\langle a,y\tr (bb') \rangle$.  It follows that there is a right multiplier $z$ of $\widehat A$ so that $y\tr b=bz$ for all $b\in\widehat A$. Using similar arguments, on finds that $z$ is actually a multiplier of $\widehat A$. This defines then a homomorphism from $\varepsilon_s(A)$ to $M(\widehat A)$. It is easy to argue that it is injective and that it maps into the multiplier algebra $\widehat A_t$ of the target algebra $\varepsilon_t(\widehat A)$. However, there seems no direct way along these lines to get that it actually has an image in the target algebra itself.
\snl
This result is known for weak Hopf algebras (see e.g.\ Lemma 2.6 in \cite{B-N-S}). 
\snl
We also get explicit formulas for these isomorphisms $\gamma_s$ and $\gamma_t$ in the case of the weak multiplier Hopf algebra associated with a separability idempotent as we study in Section \ref{s:examples}. See also Proposition 3.19 in \cite{VD8}. 
\nl
\nl

\section{\hspace{-17pt}. Special cases and examples}\label{s:examples}  

In this section, we will give some examples and use these to illustrate the results on integrals obtained in Section \ref{s:integrals} and on duality as proven in Section \ref{s:duality}. Together with these examples, we also discuss some special cases.
\snl
We begin with the two weak multiplier Hopf algebras associated with a groupoid.
\nl
\bf Integrals on $K(G)$ for a groupoid $G$ \rm
\nl
We consider any groupoid $G$. For the algebra $A$ we take the complex functions with finite support in $G$ and pointwise product. 
A coproduct $\Delta$ is defined on $A$ by 
\begin{equation*}
\Delta(f)(p,q)=
\begin{cases}
		f(pq) & \text{if $pq$ is defined},\\
				0 & \text{otherwise}.
\end{cases}
\end{equation*}
This is the first weak multiplier Hopf algebra associated with the groupoid. It is Example 1.15 in \cite{VD-W4}.
 Recall that $A_s$ is the algebra of complex functions $f$ on $G$ with the property that $f(p)=f(q)$ whenever $p,q$ are elements in $G$ with the same source. Similarly $A_t$ is the algebra of complex functions $f$ on $G$ where $f(p)=f(q)$ for elements $p,q$ with the same source. See e.g.\ a remark in the introduction of \cite{VD-W5}.
\snl
For this case, we have the following characterization of left and right integrals. Recall that for a left integral $\varphi$, we need $(\iota\ot\varphi)\Delta(f)\in A_t$ for all $f$ and that for a right integral $\psi$ we need $(\psi\ot\iota)\Delta(f)\in A_s$ for all $f$.

\stel
For any $g\in A_s$, the map $f\mapsto \sum_u g(u)f(u)$ is a left integral and any left integral on $A$ is of that form. Similarly, for any $h\in A_t$, the map $f\mapsto \sum_v h(v)f(v)$ is a right integral and any right integral on $A$ is of that form.
\estel

\bew
i) Assume that $g$ is a function on $G$ satisfying $g(p)=g(q)$ if $s(p)=s(q)$. Define $\varphi$ on $A$ by $\varphi(f)=\sum_u g(u)f(u)$. Remark that this can be defined as we assume that elements $f\in A$  have finite support. Then
\begin{equation*}
(\iota\ot\varphi)\Delta(f)(p)=\sum_u g(u)f(pu)=\sum_u g(pu)f(pu)=\sum_{u'} g(u')f(u')
\end{equation*}
where in the first two sums, we only take the sum over elements $u$ with the property that $t(u)=s(p)$ while in the last sum, we only take elements $u'$ with the property that $t(u')=t(p)$. We see immediately that this is only dependent on the target $t(p)$ and hence $(\iota\ot\varphi)\Delta(f)\in A_t$. 
\vskip 3pt
ii) Conversely, assume that $\varphi$ is a left integral. As it is a linear functional on $A$, it is of the form $f\mapsto \sum_u g(u)f(u)$ for a function $g$ on $G$. Take for $f$ the function that is $1$ in an element $q$ and $0$ in all other elements. Then 
\begin{equation*}
(\iota\ot\varphi)\Delta(f)(p)=
\begin{cases}
g(p^{-1}q) & \text{if } t(p)=t(q), \\
0 & \text{otherwise.}
\end{cases}
\end{equation*}
If this function of $p$ is only dependent on $t(p)$, that is the source of $p^{-1}$, we must have that $g(u)$ is only dependent on $s(u)$.
\vskip 3pt
iii) Similarly, for a function $h$ in $A_t$, the map $\psi$ on $A$ defined by $\psi(f)=\sum_v h(v)f(v)$ will be a right integral and any right integral is of this form.
\ebew

Observe that we can take $g=1$ and $h=1$. We will then get a faithful left integral that is also a right integral. 
\snl
On the other hand, there exists left integrals that are not right invariant and right integrals that are not left invariant. The relation between left integrals and right integrals is found in Proposition \ref{prop:2.5}. As in this example the underlying algebra is abelian, there is only one case we have to consider. We illustrate this, and some other properties, in the following example. We are still working with $A=K(G)$ as above.

\voorb i) Take a left integral $\varphi$ and a right integral $\psi$. Then $\varphi$ has the form $f\mapsto \sum_u g(u)f(u)$ where $g$ is a function in $A_s$ and $\psi$ is of the form $g\mapsto \sum_v h(v)f(v)$ with $h\in A_t$. Now take any function $c$ on $G$ and let $a=gc$ and $b=hc$. Then we have obviously
\begin{equation*}
\psi(fa)=\sum_v h(v)(fa)(v)=\sum_v h(v)f(v)a(v)=\sum_v h(v)g(v)c(v)f(v)
\end{equation*}
for all $f$. A similar argument gives the same result for $\varphi(fb)$ so that $\psi(fa)=\varphi(fb)$ for all $f\in A$. The reader can verify that the expressions for $a$ and $b$ in i) of Proposition \ref{prop:2.5}, in terms of the elements $p,q$, will yield precisely the above choice with $c$ given by $v\mapsto \sum_r p(vr^{-1})q(r)$ (and where the sum is taken over the elements $r$ with the same source as the element $v$ so that $vr^{-1}$ is defined. 
\vskip 3pt
ii) Remark that the more natural choice for the function $c$ above would be $1$. This will give the strongest relation of the form $\psi(fg)=\varphi(fh)$ for all $f\in A$. 
\vskip 3pt
iii) However, in the case e.g. that $\varphi$ is faithful, that is when $g$ is never $0$, we can take $c=g^{-1}$ so that $a=1$ and we find that $\psi(f)=\varphi(f\delta)$ where $\delta=hg^{-1}$.  Now we have $\Delta(\delta)=E(h\ot g^{-1})$ where $E$ is the canonical idempotent. This illustrates Proposition \ref{prop:2.8}.
\evoorb

Similar observations can be made for the relation between two left integrals and two right integrals. If e.g.\ $\varphi_1$ and $\varphi_2$ are left integrals, given by functions $g_1$ and $g_2$ in $A_s$ respectively, then $\varphi_1(fa)=\varphi_2(fb)$ if $a=g_2c$ and $b=g_1c$ for some $c$. The strongest choice is again with $c=1$. If however, say $\varphi_1$ is faithful, so that $h_1$ is invertible, we can take $c=g_1^{-1}$ so that $b=1$ and we get $\varphi_2(f)=\varphi_1(fa)$ with $a=g_2g_1^{-1}$. In this case, still $a\in A_t$. This illustrates Proposition \ref{prop:2.7}.
\snl
We can also illustrate the other results. Consider e.g.\ Proposition \ref{prop:2.3}. If $\varphi$ is a left integral, given by the function $g$, then for $f\in A$ and $p\in G$ we find
\begin{equation*}
((\iota\ot\varphi)\Delta(f))(p)=\sum_u g(u)f(pu)=\sum_u g(p^{-1}v)f(pp^{-1}v)=\sum_u g(v)f(pp^{-1}v)
\end{equation*}
where the first sum is taken over elements $u$ with $t(u)=s(p)$ and where the following sums taken over the elements $v$ with $t(v)=t(p)$. This last expression is equal to 
$$((\iota\ot\varphi)\Delta(f))(pp^{-1}).$$
Next, let us illustrate the formulas in Proposition \ref{prop:2.4}. Again take  a function $g\in A_s$ and the left integral $\varphi$ associated with it. Then for every $f_1,f_2\in A$ and $p\in G$ we find 
\begin{align*}
((\iota\ot\varphi)((1\ot f_1)\Delta(f_2)))(p)
&=\sum_u g(u)f_1(u)f_2(pu)\\
&=\sum_v g(p^{-1}v)f_1(p^{-1}v)f_2(v)\\
&=\sum_v g(v)f_1(p^{-1}v)f_2(v)
\end{align*}
where the first sum is taken over elements $u$ satisfying $t(u)=s(p)$ and the other sums over element $v$ satisfying $t(v)=t(p)$. The last expression is now $((\iota\ot\varphi)(\Delta(f_1)(1\ot f_2)))(p^{-1})$.
\snl
Finally, let us illustrate some of the formulas obtained in Proposition  2.7 of \cite{K-VD} and given again in Proposition \ref{prop:2.4a} in Section \ref{s:integrals}.
\nl
Because the algebra is abelian, we have that $S^{-1}=S$ and so $F_1=F_3$ and $F_2=F_4$. Hence, we have only one formula for a left integral (cf.\ Equation (\ref{eqn:2.3})) and only one formula for a right integral (cf.\ Equation (\ref{eqn:2.4})). 
\snl
If $\varphi$ is a left integral given by the function $g\in A_s$, we have on the one hand that 
\begin{equation}
(\iota\ot\varphi)\Delta(f)(p)=\sum_u g(u)f(pu)=\sum_u g(v)f(v)\label{eqn:4.6}
\end{equation}
for all $f$ and all $p$, where the first sum is taken over elements $u$ with $t(u)=s(p)$ and the second sum over the elements $v$ with the same target as $p$. On the other hand we get
\begin{align}
(\iota\ot\varphi)(F_2(1\ot f))(p)\nonumber
&=\sum_u ((S\ot\iota)E)(p,u)f(u)g(u)\label{eqn:4.7} \\
&=\sum_v f(v)g(v)
\end{align}
where the last sum is taken over those elements $v$ that satisfy $(t(v)=t(p)$. 
We see that  the formulas in (\ref{eqn:4.6}) and (\ref{eqn:4.7} yield the same result. This illustrates Equation (\ref{eqn:2.3}) of Section \ref{s:integrals}.
\snl
In a similar way, we can check the formula for a right integral (Equation (\ref{eqn:2.3})).
\nl
\bf The dual of $K(G)$ for a groupoid $G$ \rm
\nl
We again consider  a groupoid $G$ and its associated weak multiplier Hopf algebra $(A,\Delta)$ with $A=K(G)$ as in the previous item. We will now construct the dual $\widehat A$, denoted also by $B$, and by doing this, we illustrate the general procedure. As expected, we get the groupoid algebra as given already in Example 1.16 of \cite{VD-W4}. 
\snl
Of course, it would also be possible to treat this case independently, just as we did for $K(G)$. But since we want to illustrate the construction of the dual in an very explicit way, we treat it as below.
\snl
First we get the dual algebra $B$ in the next proposition.

\prop
Take the left integral $\varphi$ on $A$ defined by $f\mapsto \sum_p f(p)$. For any $h\in K(G)$ we consider the linear functional $f\mapsto \varphi(fh)$. This defines a bijective linear map from $K(G)$ to the dual space $B$. Denote by $\lambda_p$ the linear functional on $K(G)$ given by $f\mapsto f(p)$. These elements span $B$ and the product in $B$ is given by
\begin{equation*}
\lambda_p\lambda_q=
\begin{cases}
		\lambda_{pq} & \text{if $pq$ is defined},\\
				0 & \text{otherwise}.
\end{cases}
\end{equation*}
\eprop

\bew
i) Take $h\in K(G)$. Define $\gamma(h)$ as the linear functional $f\mapsto \varphi(fh)$. Clearly $\gamma(h)$ belongs to $B$. The map $\gamma$ is injective because $\varphi$ is faithful. On the other hand, when $\varphi_1$ is another left integral, then $\varphi_1(f)=\sum_u g(u)f(u)$ for some $g\in A_t$. If now  $h_1$ is any function in $K(G)$, we get
\begin{equation*}
\varphi_1(fh_1)=\sum_u g(u)f(u)h_1(u)=\varphi(f(gh_1))
\end{equation*}
for all $f$. Because $gh_1$ is again in $K(G)$, we see that $\gamma$ is also surjective. This proves the first statement.
\vskip 3pt
ii) For any $p$ and $q$ we find that 
\begin{equation*}
\langle f,\lambda_p\lambda_q \rangle = \langle \Delta(f),\lambda_p\ot\lambda_q \rangle=\Delta(f)(p,q)
\end{equation*}
and by the definition of $\Delta$ on $A$ we see that $\langle f,\lambda_p\lambda_q \rangle=f(pq)$ when $pq$ is defined and $\langle f,\lambda_p\lambda_q \rangle=0$ otherwise. This proves that $\lambda_p\lambda_q=\lambda_{pq}$ when $pq$ is defined and $\lambda_p\lambda_q=0$ if not.
\ebew

We see that the map $\gamma$ is an isomorphism of the groupoid algebra $\mathbb C G$ to $B$. 
\snl
In what follows we will identify a function $g$ in $K(G)$ with its image $\gamma(g)$ in $B$ and so we obtain the pairing 
\begin{equation*}
(f,g)\mapsto \sum_p f(p)g(q)
\end{equation*}
where $f,g$ are functions in $K(G)$ but with $f$ as an element in $A$ and $g$ as an element in $B$.
\snl
It is not hard to verify that the groupoid algebra is idempotent, and in fact,  that it has local units. Indeed, given any element $b:=\sum_p g(p)\lambda_p$ in $\mathbb C$, we can consider the sums $c$ and $d$  of the elements $\lambda_e$ and $\lambda_f$ where respectively, the first sum is taken of the sources of the elements in the support of $g$ and the second one with the targets of these elements. Then $bc=c$ and $db=b$. 
\snl
In particular, the algebra is non-degenerate. It has a unit if and only if the set of units in the groupoid is finite.
\nl
In the next proposition, we obtain the coproduct on the dual $B$ using the general procedure described in the previous section.

\prop\label{prop:4.4}
The coproduct on $B$, dual to the product on $A$, as obtained in general in Proposition \ref{prop:3.7} satisfies and is characterized by $\Delta(\lambda_p)=\lambda_p\ot\lambda_p$ for all $p\in G$.
\eprop

\bew
Let us look at the map $T_1$ on $B\ot B$.  For all $f\in K(G\ot G)$ we find
\begin{equation*}
\langle f, T_1(\lambda_p \ot \lambda_q)\rangle=T_2(f)(p,q)=f(p,pq)
\end{equation*}
if $pq$ is defined. Otherwise, we get $0$. Hence $T_1(\lambda_p \ot \lambda_q)=\lambda_p \ot \lambda_{pq}$ if $pq$ is defined and we get $0$ otherwise. It follows that $\Delta(\lambda_p)=\lambda_p \ot \lambda_p$ for all $p$.
\ebew

It is immediately clear from this property that $\Delta$ is a homomorphism and that it is coassociative. Also fullness of $\Delta$ is a consequence. In this case, the coproduct is coabelian. 
For the antipode, we obviously have $S(\lambda_p)=\lambda_{p^{-1}}$. The square of the antipode is again the identity map. 
For the source and target maps and the source and target algebras, we have the following result.

\prop
The source and target maps on the dual $B$ are given by
\begin{equation*}
s(\lambda_p)=\lambda_{s(p)}
\tussen 
t(\lambda_p)=\lambda_{t(p)}
\end{equation*}
for all $p$. The source and target algebras $\varepsilon_s(B)$ and $\varepsilon_t(B)$ are the same and spanned by elements $\lambda_e$ where $e$ is a unit. There multiplier algebras $B_s$ and $B_t$ are identified with all complex functions on the units, with pointwise product.
\eprop

We can now illustrate Proposition \ref{prop:3.19}. 
\snl
For all $y\in \varepsilon_s(A)$ we should have an element $\gamma_s(y)$ in the target algebra of $B$ satisfying $\langle ya,b\rangle=\langle a,b\gamma_s(y)\rangle$  for all $a\in A$ and $b\in B$. Take any element $a\in A$ while  $b=\lambda_p$ where $p\in G$. For the left hand side we get $a(p)y(p)$ and for the right hand side $\sum_r a(r)(b\gamma_s(y))(r)$. This must hold for all $a$ and it follows that 
$(b\gamma_s(y))(p)=y(p)$ and that it is $0$ in other points. Because $b=\lambda_p$, this means that $\gamma_s(y)$ must have support in the units and that in fact $\gamma_s(y)(e)=y(r)$ if $e=s(r)$. This is well-defined because it is assumed that $y(r)=y(r')$ if $r$ and $r'$ have the same source. One can verify that the map $\gamma_s$ is indeed an isomorphism from the source algebra  of $A$ to the target algebra of $B$.
\snl
Similarly the isomorphism $\gamma_t$ of the target algebra of $A$ to the source algebra of $B$ is given here by $\gamma_t(x)(e)=x(r)$ when $e=t(r)$.
\snl
The general formula (cf.\ \ref{prop:3.9}) for the dual counit yields that $\varepsilon(\lambda_p)=1$ for all $p$ in $G$. This is indeed the unique counit on the dual $(B,\Delta)$. In the following proposition, we obtain the dual canonical idempotent $\widehat E$.

\prop
The canonical idempotent for the dual is $\sum_e \lambda_e\ot \lambda_e$ where the sum is taken over the units of $G$.
\eprop

\bew
We have seen in \ref{prop:3.13} that $\widehat E$ is given by the formula $\langle a\ot a',\widehat E\rangle=\varepsilon(aa')$ where $\varepsilon$ is the counit of $A$. The counit on $K(G)$ is obtained by taking the sum of the values in the units. So here we get the sum  $\sum_e a(e)a'(e)$ over the units $e$. This means that, as a function of two variables, we have $\widehat E(e,e)=1$ and $0$ in all other points. Hence $\widehat E=\sum_e \lambda_e\ot \lambda_e$.
\ebew

Of course, this is what we expect as the identity $1$ in $M(B)$ is the sum $\sum_e\lambda_e$ and so $\Delta(1)=\sum_e \lambda_e\ot\lambda_e$.
\snl
It is instructive to verify that this sum indeed defines a multiplier of $B\ot B$ and that it is the smallest idempotent in $M(B\ot B)$ with the property that $\widehat E\Delta(b)=\Delta(b)$ and $\Delta(b)\widehat E=\Delta(b)$. For this, it is enough to look at elements $b$ of the form $\lambda_p$ with $p\in G$. Similarly, one can verify that $\widehat E(B\ot B)=\Delta(B)(1\ot B)=\Delta(B)(B\ot 1)$ as well as $(B\ot B)\widehat E=(1\ot B)\Delta(B)=(B\ot 1)\Delta(B)$.

\nl
\bf Integrals on the dual \rm
\nl
It is not hard to find the integrals on the dual $B$ directly from the definition of the coproduct $\Delta$ on $B$ as obtained in Proposition \ref{prop:4.4}. However, following the spirit of this section, we will obtain the integrals from the general result we found in the section on duality, in Proposition \ref{prop:3.16}. Remark that left integrals are also right integrals and vice versa because the coproduct is coabelian.

\prop
If $g$ is a function on $G$ with support in the set of units, then $\lambda_p\mapsto g(p)$ is an integral on $B$. Any integral is of this form.
\eprop

\bew
i) As in Proposition \ref{prop:3.16} take $a\in A$ and define $\psi_a$ on $B$ by
\begin{equation*}
	\psi_a(\omega)=\varphi(a\varepsilon_s(c))\qquad\text{ when } \qquad \omega=\varphi(\,\cdot\, c) 
\end{equation*}
where $\varphi$ is a left integral on $A$ and $c$ an element $A$. We can take for $\varphi$ the integral $f\mapsto \sum_p f(p)$. Given the identification of functions in $K(G)$ with elements in $B$, we get
$\psi_a(f)=\sum_u a(u)f(s(u))$, in particular $\psi_a(\lambda_p)=\sum_u a(u)$ where the sum is taken over the elements $u$ with $s(u)=p$.  We see that $\psi_a(\lambda_p)$ is $0$ except if $p$ is a unit. 
\vskip 3pt
ii) Conversely, for any function $g$ on $G$ with support in the set of units, the map $\psi:\lambda_p\mapsto g(p)$ will be an integral. Indeed, for such a function we have 
\begin{equation*}
(\psi\ot\iota)\Delta(\lambda_p)=g(p)\lambda_p
\end{equation*}
and this belongs to $B_s$ as this is precisely the algebra of functions with support in the units.
\ebew
Also here, we have a faithful integral:

\prop Define $\varphi$ on $B$ by $\varphi(\lambda_e)=1$ when $e$ is a unit and $\varphi(\lambda_p)=0$ otherwise. Then $\varphi$ is a faithful integral.
\eprop

\bew
We know already that $\varphi$, defined as in the formulation of the proposition, is an integral. We just have to show that this one is faithful.
\snl
So, let $b$ be any element in $B$ and assume that $\varphi(\lambda_p b)=0$ for all $p$. Now, given $p$, we have that  $pq$ is a unit if and only if $q=p^{-1}$. It follows that $b(p^{-1})=0$. This holds for all $p$ and so $b=0$. Similarly on the other side.
\ebew

This is the integral that we get by applying the result given of Proposition \ref{prop:2.19}. And it illustrates Theorem \ref{stel:2.20}.

\nl
\bf The weak multiplier Hopf algebra associated with a separability idempotent \rm
\nl
In this item we briefly treat a special example of a weak multiplier Hopf algebra associated with a separability idempotent. In particular, the dual weak multiplier Hopf algebra for this example is constructed. However, we do not give details here. We refer to a separate note  where this example and its dual are treated with full details (see \cite{VD8}). 
\snl
See also Item 5.3 in Section 5 of \cite{T2} for a related example.
\nl
The starting point is a pair $B, C$ of non-degenerate and idempotent algebras and a separability idempotent $E\in M(B\ot C)$. Recall  Definition 1.4 from \cite{VD4}.

\defin 
Let $E$ be an idempotent in the multiplier algebra $M(B\ot C)$ and assume that $E(1\ot b)$ and $(c\ot 1)E$ belong to $B\ot C$ for all $b\in B$ and $c\in C$. We assume that $E$ is full in the sense that the left leg and the right leg of $E$ are respectively  all of $B$ and all of $C$. Furthermore it is required that there are non-degenerate anti-homomorphisms $S_B:B\to M(C)$ and $S_C:C\to M(B)$ satisfying (and characterized by)
\begin{equation*}
E(b\ot 1)=E(1\ot S_B(b))
\tussenen
(1\ot c)E=(S_C(c)\ot 1)E
\end{equation*}
for all $b\in B$ and $c\in C$. Then $E$ is called a separability idempotent. 
\edefin

If the maps $S_B$ and $S_C$ have range in $C$ and $B$ respectively, $E$ is called semi-regular and if moreover these maps are anti-isomorphisms, then $E$ is called regular. Here we will  assume that $E$ is regular. 
\snl
The regular case was first studied in the first version of \cite{VD4}, while the more general case is considered in the second version of this paper.
\snl
We denote by $\varphi_B$ and $\varphi_C$ the distinguished linear functionals on $B$ and $C$  satisfying
\begin{equation*}
(\varphi_B\ot \iota)E=1
\tussenen
(\iota\ot\varphi_C)E=1
\end{equation*}
in $M(C)$ and $M(B)$ respectively. Recall that in the regular case, the distinguished functionals are faithful and have KMS automorphisms (denoted by $\sigma_B$ and $\sigma_C$ respectively), see Proposition 1.13 in \cite{VD4}. Recall that $\sigma_B$ is the inverse of $S_CS_B$ while $\sigma_C=S_BS_C$.
\snl
We have the following associated weak multiplier Hopf algebra (cf.\ Proposition 3.2 in \cite{VD-W5}). 

\prop Let $A$ be the algebra $C\ot B$. Define $\Delta$ on $A$ by
\begin{equation*}
\Delta(c\ot b)=c\ot E\ot b
\end{equation*}
for all $b\in B$ and $c\in C$. Then we have a weak multiplier Hopf algebra. 
\vskip 3pt
The counit is given by $\varepsilon(c\ot b)=\varphi_B(S_C(c)b)=\varphi_C(cS_B(b))$. The antipode satisfies $S(c\ot b)=S_B(b)\ot S_C(c)$ and the source and target maps are 
\begin{equation*}
\varepsilon_s(c\ot b)=1\ot S_C(c)b
\tussenen
\varepsilon_t(c\ot b)=cS_B(b)\ot 1
\end{equation*}
for all $b,c$.
\eprop

We see that $\varepsilon_s(A)=1\ot B$ and $\varepsilon_t(A)=C\ot 1$. The multiplier algebras $A_s$ and $A_t$ are $1\ot M(B)$ and $M(C)\ot 1$ respectively. 
\snl
Observe also that $A$ is regular if and only if $E$ is regular (in the sense of Definition 2.3 of \cite{VD4}. See Proposition 3.3 in \cite{VD-W5}.
\snl
In what follows we will look at $A$ as the algebra generated by $B$ and $C$ with the condition that elements of $B$ commute with elements of $C$. We can do this by looking at the natural non-degenerate homomorphisms from $B$ and $C$ into $M(A)$ and by identifying $B$ and $C$ with their images. Then $A$ is spanned by elements of the form $cb$ and $cb=bc$ in $A$. This practice is standard. 
\snl
Then we write e.g.\ $\Delta(cb)=(c\ot 1)E(1\ot b)$ and the canonical idempotent of $A$ is $E$, considered as sitting in $M(A\ot A)$. We refer to Section 3 in  \cite{VD-W5}.
\snl
We get the following characterization of the integrals. See Proposition 2.2 in \cite{VD8}.

\prop
For any linear functional $g$ on $B$ we have a left integral $\varphi$ on $A$ given by $\varphi(cb)=\varphi_C(c)g(b)$ and  any left integral has this form. Similarly, for any linear functional $f$ on $C$ we have a right integral $\psi$ on $A$ given by $\psi(cb)=f(c)\varphi_B(b)$ and again any right integral is of this form. Moreover, the map $cb\mapsto \varphi_C(c)\varphi_B(b)$ is a faithful left integral that is also a right integral.
\eprop

The proof of this result is straightforward, again see \cite{VD8}. 
\nl
\bf The dual of $(CB,\Delta)$ \rm
\nl 
First observe that we use $B$ and $C$ for the underlying algebras. So we will not use $B$ for the dual $\widehat A$ here in this item.
\snl
Because we have a faithful integral $\varphi$ on $A$, defined by $\varphi(cb)=\varphi_C(c)\varphi_B(b)$, we have that $\widehat A$, as a vector space, can be identified with $C\ot B$ via the map 
$v\ot u\mapsto \varphi(\,\cdot\,vu)$. Remark that $\varphi(cbvu)=\varphi_C(cv)\varphi_B(bu)$ for all $b,u\in B$ and $c,v\in C$. So this map is $v\ot u\mapsto \varphi_C(\,\cdot\,v)\ot \varphi_B(\,\cdot\,u)$. 
\snl
However, in what follows, we will make another identification. It turns out to give nicer formulas.

\defin\label{defin:pairing}
We define a pairing of vector spaces $A$ and $B\ot C$ by
\begin{equation*}
\langle cb,u\ot v \rangle = \varphi_B(bS_C(v))\varphi_C(S_B(u)c)).
\end{equation*}
\edefin

Observe the difference with the formula for the counit on $A$. Indeed, for the counit we have $\varepsilon(cb)=\varphi_B(S_C(c)b)=\varphi_C(cS_B(b))$.
\snl
Because we have
\begin{align*}
 \varphi_B(bS_C(v))\varphi_C(S_B(u)c)
 &=\varphi_B(bS_C(v))\varphi_C(c\sigma_C(S_B(u))) \\
 &=\varphi(cbS_C(v)\sigma_C(S_B(u))),
\end{align*}
we see that 
\begin{equation*}
u\ot v\mapsto \varphi(\,\cdot\,S_C(v)\sigma_C(S_B(u)))
\end{equation*}
gives a linear isomorphism from the space $B\ot C$ to $\widehat A$.
\snl
In what follows, we will make this identification, but we will use $u\diamond v$ for $u\ot v$ when we consider it as an element in $\widehat A$. We will also systematically use the letters $u,u',\dots$ and $v,v'\dots$ for elements in $B$ and $C$ respectively, when $B\ot C$ is identified with $\widehat A$ via the pairing in Definition \ref{defin:pairing} above.
\snl
For the product on the dual $\widehat A$ we find the following result. See Proposition 3.2 in \cite{VD8}.

\prop\label{prop:4.10} 
Let $u,u'\in B$ and $v,v'\in C$. For the product in $\widehat A$ we find
\begin{equation*}
(v\diamond u)(v'\diamond u')=\varepsilon(v'u)\,\,v\diamond u',
\end{equation*}
were $\varepsilon$ is the counit on $A$
\eprop

Again the proof is straightforward and it is found in \cite{VD8}.
\nl
By the general theory, we know that this product is associative and non-degenerate. Here, we can verify that the product is associative by a simple calculation. We see also that it is non-degenerate  because the map $(u,v)\mapsto \varepsilon(vu)$ is a non-degenerate pairing of the space $B$ with $C$. This follows from the definition $\varepsilon(vu)=\varphi_B(S_C(v)u)$, the fact that $\varphi_B$ is faithful and that $S_C$ is bijective.
\snl
In what follows, we will also write $B\di C$ for the algebra we obtain in Proposition \ref{prop:4.10}. The algebra is a  (possibly) {\it infinite matrix algebra}, build with two vector spaces and a non-degenerate pairing. The algebra structure is not dependent on the multiplications in $B$ and $C$. It only depends on the pairing of the underlying vector spaces.
\snl
Because the algebra $B\di C$ is non-degenerate, we can consider  its multiplier algebra \\
$M(B\di C)$. We get the following characterization. See Proposition 3.3 in \cite{VD8}. 

\prop
A linear map $\gamma:B\to B$ is called {\it adjointable} if there is a linear mapping $\gamma^t:C\to C$ satisfying $\varepsilon(v\gamma(u)=\varepsilon(\gamma^t(v)u)$ for all $u\in B$ and $v\in C$. The maps $\gamma$ and $\gamma^t$ determine each other. For any adjointable map $\gamma:B\to B$ with adjoint $\gamma^t$ there is a multiplier $m$ of $B\di C$ given by
\begin{equation}
m(u\di v)=\gamma(u)\di v
\tussenen
(u\di v)m=u\di\gamma^t(v). \label{eqn:multiplier}
\end{equation}
Any multiplier is of this form.
\eprop

The embedding of $B\di C$ in $M(B\di C)$ is found by associating the linear map $u'\mapsto \varepsilon(vu')u$ to the element $u\di v$. It adjoint is $v'\mapsto \varepsilon(v'u)v$. The identity map from $B$ to itself is of course adjointable and the associated multiplier is the identity $1$ in $M(B\di C)$
 \nl
 \bf The coproduct on $B\di C$\rm
 \nl
 By definition, the coproduct on $\widehat A$ is dual to the product on $A$. Because $A$ is the tensor product of the algebras $C$ and $B$, it is expected that the coproduct on the dual is also a tensor product of coproducts on the factors. In this case, it would mean that the coproduct on $B\di C$ has the form 
 \begin{equation*}
\Delta(u\di v)=\sum_{(u),(v)} (u_{(1)}\di v_{(1)})\ot (u_{(2)}\di v_{(2)})
\end{equation*}
where we have the Sweedler notation for coproducts $\Delta_B$ on $B$ and $\Delta_C$ on $C$.
\snl
This is indeed the case and it is made precise with the following results.

\prop\label{prop:4.14} 
Define a pairing of $C$ with $B$ by the formula $\langle c,u \rangle_1=\varphi_C(S_B(u)c)$. For this pairing we have
\begin{equation*}
\langle cc',u \rangle_1 = \langle c\ot c',\Delta_B(u) \rangle_1
\end{equation*}
for $c,c'\in C$ and $u\in B$, where $\Delta_B(u)=F_1(1\ot u)$ and $F_1=(\iota\ot S_C)E$ as in Section \ref{s:integrals}.
\eprop

\prop\label{prop:4.15}
Define a pairing of $B$ with $C$ by the formula $\langle b,v \rangle_2= \varphi_B(bS_C(v))$. For this pairing we have 
\begin{equation*}
\langle bb',v\rangle_2 =\langle b\ot b',\Delta_C(v)\rangle_2
\end{equation*}
where $\Delta_C(v)=(v\ot 1)F_2$ for $b,b'\in B$ and $v\in C$ and $F_2=(S_B\ot\iota)E$ as in Section \ref{s:integrals}.
\eprop
 
Now we get from these two properties the formula for the coproduct on $B\di C$.

\prop \label{prop:4.16}
The coproduct on $B\di C$, defined by the formula 
\begin{equation*}
\langle aa',u\di v\rangle=\langle a\ot a',\Delta(u\di v)\rangle
\end{equation*}
 for $a,a'\in A$ and $u\in B$ and $v\in C$, is
 \begin{equation*}
\Delta(u\di v)=\sum_{(u),(v)} (u_{(1)}\di v_{(1)})\ot (u_{(2)}\di v_{(2)})
\end{equation*}
where we use the Sweedler notations
\begin{align*}
\sum_{(u)} u_{(1)}\ot u_{(2)}&=\Delta_B(u)=F_1(1\ot u)\\
\sum_{(v)} v_{(1)}\ot v_{(2)}&=\Delta_C(v)=(v\ot 1)F_2.
\end{align*}
\eprop 

We need to make a couple of remarks about this formula.

\opm
i) It is rather remarkable that the pairings $(c,u)\mapsto \langle c , u\rangle_1$ and $(b,v)\mapsto \langle b, v\rangle_2$, that we defined in the Propositions \ref{prop:4.14} and \ref{prop:4.15}, yield coproducts $\Delta_B:B\to B\ot B$ and $\Delta_C:C\to C\ot C$. 
\vskip 2pt
ii) Moreover, these coproducts are known. See e.g.\ Proposition 2.9 in  \cite{VD4}.
\vskip 2pt
iii) It is also obvious that the coproducts are homomorphisms. As a consequence, the tensor coproduct, as defined in Proposition \ref{prop:4.16} will be a homomorphism on the tensor product algebra $B\ot C$. However, this is different from the product in $B\di C$ and, from the general theory, we know that this tensor coproduct is also a homomorphism on $B\di C$.  
\eopm

\snl
\bf The regular weak multiplier Hopf algebra $B\di C$ \rm
\nl
One can now systematically verify that indeed, the algebra $B\di C$ as obtained in Proposition \ref{prop:4.10}, together with the coproduct given in Proposition 4.16, satisfies all the requirements of a regular weak multiplier Hopf algebra. This is done in detail in Section 3 of \cite{VD8}. 
\snl
For the counit  $\widehat\varepsilon$ on the dual we find that $\widehat\varepsilon(u\di v)=\varphi_B(u)\varphi_C(v)$. This is no surprise because $\varphi_B$ and $\varphi_C$ are the counits for the coalgebras $(B,\Delta_B)$ and $(C,\Delta_C)$.
\snl
For the antipode $\widehat S$ on the dual we also get the expected formula. We have
\begin{equation*}
\widehat S(u\di v)=S_B^{-1}(v)\di S_C^{-1}(u)
\end{equation*}
for $u\in B$ and $v\in C$. This result is proven by taking the adjoint of the antipode in the pairing. See Proposition 3.10 in \cite{VD8}. 
\snl
One can verify that $\widehat S$ is an anti-isomorphism of the algebra $\widehat A$ and that it flips the coproduct. See Propositions 3.11 and 3.12 in \cite{VD8}. 
\snl
The canonical idempotent is more difficult to handle because it is an element in the multiplier algebra $M(B\di C)$ of $B\di C$. There is the following result (see Proposition 3.13 in \cite{VD8}).

\prop\label{prop:4.18}
The maps $\gamma: u\ot u'\mapsto (uu'\ot 1) F_1$ from $B\ot B$ to itself and  $\gamma':v\ot v'\mapsto F_2(1\ot vv')$ from $C\ot C$ to itself are adjoint to each other (with respect to the tensor product pairing). They are projection maps. And they give the canonical idempotent $\widehat E$ in $M(\widehat A\ot\widehat A)$.
\eprop

That the maps are idempotent follows simply from the fact that $m_BF_1=E_{(1)}S_C(E_{(2)})=1$ as well as $m_CF_2=S_B(E_{(1)})E_{(2)}$. Here we use $m_B$ and $m_C$ for the multiplication maps on $B\ot B$ and $C\ot C$ respectively.
\snl

Using the formula for the antipode, we find the following formulas for the source and target maps on the dual. 

\prop
For all $u\in B$ and $v\in C$ we have
\begin{equation*}
\varepsilon_s(u\di v)=\varphi_B(u)(E_{(1)}\di E_{(2)}v)
\tussenen
\varepsilon_t(u\di v)=\varphi_C(v)(uE_{(1)}\di E_{(2)})
\end{equation*}
where we again use the Sweedler type notation for $E$. 
\eprop

Recall that $(u\ot 1)E\in B\ot C$ when $u\in B$ and this gives the element $uE_{(1)}\di E_{(2)}$ in $B\di C$. Similarly $E(1\ot v)\in B\ot C$ when $v\in C$ and this gives the element $(E_{(1)}\di E_{(2)}v)$ in $B\di C$.
\snl
One can verify that the maps $v\mapsto E_{(1)}\di E_{(2)}v$ and $u\mapsto uE_{(1)}\di E_{(2)}$ are isomorphisms from the algebras $C$ and $B$ to the source algebras $\varepsilon_s(\widehat A)$ and $\varepsilon_t(\widehat A)$ respectively, see Proposition 3.20 in \cite{VD8}.
\snl
They are in fact modified forms of the isomorphisms obtained in general in Proposition \ref{prop:3.19}. Again see \cite{VD8} for the relation between the two results. 
\nl
Finally we consider the integrals on the dual $\widehat A$. From the general theory we find that there is a right integral $\widehat \psi_c$ on the dual given by $\widehat\psi_c(u\di v)=\varphi_C(S_B(u)cv)$ for all $u\in B$ and $v\in C$ (see Proposition 3.23 in \cite{VD8}). In fact, we get the following.

\prop
For any $x\in M(C)$ there is a right integral $\widehat\psi_x$ on $\widehat A$ given by 
\begin{equation*}
\widehat\psi_x(u\di v)=\varphi_C(S_B(u)xv)
\end{equation*}
 for all $u\in B$ and $v\in C$. Any right integral is of this form. Similarly, for any $y\in M(B)$ there is a left integral $\widehat\varphi_y$ on $\widehat A$ given by $\widehat\varphi_y(u\di v)=\varphi_B(uyS_C(v)$ for all $u\in B$ and $v\in C$. Again any left integral is of this form.
\eprop

For the proof we refer to Propositions 3.23 and 3.24 in \cite{VD8}. We also found and expression for the modular element (see Proposition 3.28 in \cite{VD8}).

\prop 
If $\widehat\varphi$ is the left integral given by $\widehat\varphi(u\di v)=\varphi_B(uS_C(v))$, then
 \begin{equation*}
\widehat\varphi(S(u\di v))=\widehat\varphi((u\di v)\delta))
\end{equation*}
 where $\delta$ is the multiplier of $\widehat A$ given by 
\begin{equation*}
\delta(u\di v)=\sigma_B^{-2}(u)\di v
\tussenen
(u\di v)\delta=u\di \sigma_C^{-2}(v)
\end{equation*}
\eprop

We finish with a remark. From the formula for the antipode we see that 
\begin{equation*}
S^2(u\di v)=S_B^{-1}S_C^{-1}(u)\di S_C^{-1}S_B^{-1}(v)=\sigma_B(u)\sigma_C^{-1}(v)
\end{equation*}
for all $u,v$. Then using the formulas for $\delta$ obtained in the previous proposition, we fin 
\begin{equation*}
S^4(u\di v)=\delta^{-1}(u\di v)\delta
\end{equation*}
for all $u,v$. This is a special form of Radford's formula as there is left integral on $A$ that is also right invariant.
\nl\nl

\section{\hspace{-17pt}. Conclusion and further research}\label{s:conclusions}  

In this paper, we have continued our work on weak multiplier Hopf algebras, started in \cite{VD-W3}, \cite{VD-W4} and \cite{VD-W5}, with the study of integrals and duality. 
\snl
Some questions remain open.
\snl
In the examples we considered in Section 4, there always existed a single faithful left integral. It is known from the theory of weak Hopf algebras, that it is possible however that there is a faithful set of integrals but not a single faithful integral. This follows  e.g.\  from Proposition 2.5 in \cite{I-K}. We are indebted to G.\ B\"ohm for informing us about this example in the paper.
\snl 
It is known that, for finite-dimensional weak Hopf algebras, a single faithful integral exists if and only if the underlying algebra is Frobenius (see Theorem 3.16 in \cite{B-N-S}). For (possibly non finite-dimensional) weak multiplier Hopf algebras, the situation is more subtle. One can expect that when the underlying algebra is Frobenius and when there are enough integrals, there will also exist a single faithful one. Moreover, given the fact that for multiplier Hopf algebras, integrals are unique, and hence the existence of enough integrals implies the existence of a single faithful one, one might wonder if this should eventually only be dependent on the base algebra?  
In fact, if we look again at the arguments at the end of Section \ref{s:integrals}, they seem to support such a conjecture.
 \snl
 In the case of a single faithful integral, we have objects like the Radon-Nikodym elements, giving any left integral in terms of a single faithful left integral (see Proposition \ref{prop:2.7}), the modular elements, expressing any right integral in terms of a single faithful integral (see Proposition \ref{prop:2.8}) and finally the modular automorphisms of a single faithful left integral (see Proposition \ref{prop:2.7a}). 
One should also prove results like the ones above in the event that there is a faithful set of integrals, but no single faithful integral.
\snl 
We have not considered the involutive case in this paper. It can be verified that the dual of a weak multiplier Hopf $^*$-algebra with integrals is again, in a natural way, a weak multiplier Hopf $^*$-algebra with integrals. It is expected that the dual will carry a faithful set of positive integrals when this is the case for the original $^*$-algebraic quantum groupoid. In that case, it should be possible to lift the whole structure to a Hilbert space setting. This is work in progress, cf.\ \cite{K-VD1} and \cite{K-VD2}. 
\snl
We also have only considered the regular case. There are reasons to believe that the existence of integrals can only be studied under the reasonable condition of regularity, but still, the problem is open. It is not entirely clear what remains possible in the non-regular case. This has not even been studied for multiplier Hopf algebras.
\snl
In our work, we only have algebras over the field of complex numbers. We believe that many results will still be valid for algebras over other fields, but this should be carefully investigated. 
\nl
Finally, as we mentioned already in the introduction, while preparing this work, Timmermann obtained results about duality for regular multiplier Hopf algebroids with integrals in \cite{T1} and \cite{T2}. And he concluded some results about integrals and duality for regular weak multiplier Hopf algebras. 
\snl
The theory of multiplier Hopf algebroids is more general, but also much more complicated than the theory of weak multiplier Hopf algebras. In order to obtain the results on duality of weak multiplier Hopf algebras with integrals from the theory of multiplier Hopf algebroids, one needs many different concepts and results. There is not only the theory of integrals and the duality for multiplier Hopf algebroids (as given in \cite{T1} and \cite{T2}). One also needs to understand the procedure to pass from a weak multiplier Hopf algebra to a multiplier Hopf algebroid and back (as in \cite{T-VD1}). And finally one needs definitions and results on integrals on weak multiplier Hopf algebras (as already available in \cite{K-VD}). 
\snl
On the other hand, the direct approach to the duality for weak multiplier Hopf algebras, as developed in this paper, is less complicated. Moreover, in this paper, we also give more results and also more details about integrals on weak multiplier Hopf algebras. And further, the example we treat here in Section \ref{s:examples} (and in greater detail in \cite{VD8}) provides a deeper insight in the theory as it is a careful direct application to illustrate our results.
\snl
Still, we believe it is interesting to do some more research about where and how results on integrals and duality for multiplier Hopf algebroids from \cite{T1} and \cite{T2} are related with more results from this paper.
\nl\nl

\end{document}